\documentclass{article}

\pdfoutput=1

\usepackage{natbib}
\usepackage[normalem]{ulem}
\usepackage{graphicx}
\usepackage{amsmath,amsthm,bm,mathrsfs}
\usepackage{amsfonts}
\usepackage{booktabs}
\usepackage[table]{xcolor}
\usepackage{pgfplots,tikz}
\usepackage{rotating}
\usepackage{multirow}
\usetikzlibrary{patterns,fadings}
\usepackage{pgfplots}
\usepackage{makecell}
\usetikzlibrary{arrows}
\usepackage[titletoc]{appendix}

\usepackage[margin=1.4in]{geometry}

\newcommand{\Frac}[2]{\displaystyle \frac{#1}{#2}}
\newcommand{\deriv}[3]{\Frac{\partial^{#1} #2}{\partial {#3}^{#1}}}

\newcommand{\nabu}{{\boldsymbol \nabla}}
\newcommand{\grad}{\nabu}
\newcommand{\Div}{\grad \cdot}

\renewcommand{\vector}[1]{\underline{\boldsymbol{#1}}}

\newcommand{\dint}{\displaystyle\int}

\definecolor{green}{rgb}{0.034, 0.5, 0.034}

\newcommand{\rev}[1]{{\textcolor{black}{#1}}}
\newcommand{\revtwo}[1]{{\textcolor{black}{#1}}}

\renewcommand{\theequation}{\thesection.\arabic{equation}}
\numberwithin{equation}{section}

\pgfmathsetmacro{\radius}{0.2}
\pgfmathsetmacro{\length}{8}

\begin{document}

%%%%%%%%%%%%%%%%%
\title{Emergent three-dimensional sperm motility: \\Coupling calcium dynamics and preferred curvature \\ in a Kirchhoff rod model}
\author{ {\sc Lucia Carichino and Sarah D. Olson}\\[2pt]
Department of Mathematical Sciences, Worcester Polytechnic Institute, \\
100 Institute Rd, Worcester MA 01609, USA.\\[6pt]
{\rm Published: 16 October 2018}\vspace*{6pt}\\[6pt]
https://doi.org/10.1093/imammb/dqy015
}
\pagestyle{headings}
\markboth{L. CARICHINO}{\rm 3D Sperm Motility: Calcium-Curvature Coupling}
\maketitle

%%%%%%%%%%%%%%%%%abstract style
\begin{abstract}
{Changes in calcium concentration along the sperm flagellum \revtwo{regulate} sperm motility and hyperactivation, characterized by an increased flagellar bend amplitude and beat asymmetry, enabling the sperm to reach and penetrate the ovum (egg). The signaling pathways by which calcium increases within the flagellum are well established. However, the exact mechanisms of how calcium regulates flagellar bending are still under investigation. We extend our previous model of planar flagellar bending by developing a fluid-structure interaction model that couples the three-dimensional motion of the flagellum in a viscous, Newtonian fluid with the evolving calcium concentration. The flagellum is modeled as a Kirchhoff rod: an elastic rod with preferred curvature and twist. The calcium dynamics are represented as a one-dimensional reaction-diffusion model on a moving domain, the centerline of the flagellum. The two models are coupled assuming that the preferred curvature and twist of the sperm flagellum depend on the local calcium concentration. To investigate the effect of calcium on sperm motility, we compare model results of flagellar bend amplitude and swimming speed for three cases: planar, helical (spiral with equal amplitude in both directions), and quasi-planar (spiral with small amplitude in one direction). 
We observe that for the same parameters, the planar swimmer is faster and a turning motion is more clearly observed when calcium coupling is accounted for in the model. In the case of flagellar bending coupled to the calcium concentration, we observe emergent trajectories that can be characterized as a hypotrochoid for both quasi-planar and helical bending.}
{Sperm motility; Kirchhoff rod; Calcium dynamics; Hypotrochoid; Regularized Stokeslets}
\end{abstract}
%%%%%%%%%%%%%%%%%%%%%%%%%%%%%%%
%%%%%%%%%%%%%%section A%%%%%%%%%
\section{Introduction}\label{sec:intro}
Changes in cytosolic calcium concentration in the sperm flagellum have been shown to be essential to achieve fertilization of the ovum (egg) ~\citep{Guerrero10,Guerrero11,Suarez08}. An increase in calcium concentration is associated with hyperactivated motility~\citep{HoSuarez01a,Suarez03}, which is characterized by an increased flagellar bend amplitude and beat asymmetry. This motility pattern  enables the sperm to escape from the oviductal sperm reservoir~\citep{Demott_92}, to detach from the oviductal epithelium~\citep{Demott_92}, to switch from a straight motion to a curved motion and back~\citep{Darson_08}, and to penetrate the egg~\citep{Drobnis88,Quill_03,Ren_01}.

In particular, in mammalian sperm, the opening of CatSper channels on the principal piece of the flagellum has been associated with an increase in intracellular calcium concentration and hyperactivated motility~\citep{Carlson03,HoSuarez01b,Ho02,Suarez03,xia_07}. Experiments have shown that CatSper-null mutant sperm \revtwo{have impaired motility, which correlates with} infertility~\citep{Ho09}. 
The exact mechanisms through which the increase in calcium concentration influences flagellar bending at the level of dyneins, the active force generators, is not completely known. Calcium is hypothesized to bind to centrin or calmodulin, different calcium binding proteins within the axonemal structure of the flagellum \citep{Lindemann95}.  The increase in calcium does change the bending, but the exact timing and magnitude of force generation is not completely known. However, in a phenomenological model, we can explore different ways to couple evolving calcium dynamics to flagellar bending.

\begin{figure}[!tbh] %F1
\centering
(a)\hspace{2.65in}(b)
\includegraphics[width=0.8\linewidth]{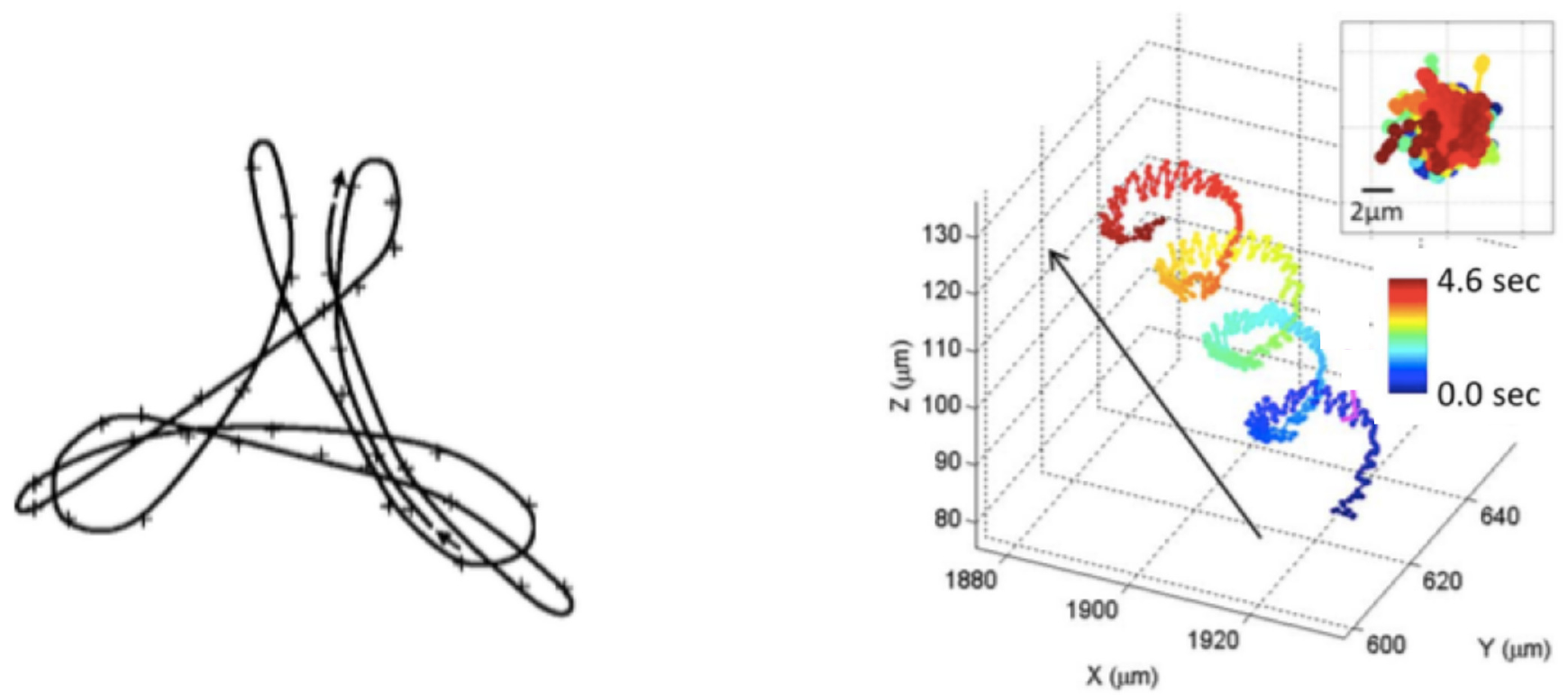}
\caption{\textit{Experimental trajectories.} (a) Trace of the mid-piece of a mouse sperm swimming near a cover slip for a time period corresponding to 3.5 beats (approximately 0.35 seconds). (b) Trace of the head of a horse sperm for a time period of 4.6 seconds. The figures are reproduced, with permission, from \cite{Woolley03} and \cite{Su13}, respectively.}
\label{exppic}
\vspace*{-9pt}
\end{figure}

Sperm are navigating in a complex, three-dimensional (3D) fluid environment in order to reach and penetrate the egg. Thus, it is necessary to study sperm motility in 3D ~\citep{Guerrero11}. However, many experiments have analyzed the motility of a few sperm where motility is recorded at a fixed depth. In this case, motility patterns observed were restricted in a given two-dimensional (2D) plane. Recently, new technologies have been developed to trace the 3D trajectories of multiple sperm at the same time~\citep{Jikeli_15,Su12,Su13}. Observed sperm trajectories can vary form planar to quasi-planar, and to helical, possibly forming a chiral ribbon or a flagelloid curve, as shown in Figure~\ref{exppic}~\citep{smith_09,Woolley01,Woolley03, Su12,Su13}. The particular trajectory observed depends on the species of sperm, on the proximity to the oviductal walls and on the external fluid properties. In particular, \cite{smith_09} and \cite{Woolley03} showed that the flagellar beat form can change from 3D to 2D as the fluid viscosity increases. 
%%%

Different 2D and 3D elastohydrodynamic models have been developed to study flagellar motility and the interaction with the surrounding fluid, where the flagellum is either described at the continuum level ~\citep{gadelha10,Ishimoto18,olson_11,olson_13,simons_14,simons_15}, or at a more detailed level including the mathematical description of the discrete dynein motors, along with the accessory structures such as the  microtubules and nexin links in the flagellum~\citep{Yang_08}. In these models, the flagellar beat form is an emergent property. Another approach is to prescribe the flagellar beat form~\citep{curtis_12,Ishimoto_16,Ishimoto_17}. In all of these modeling studies, regardless of how the flagellar beat form is modeled, the observed sperm trajectory is an emergent property of the coupled system accounting for the fluid dynamics and the swimmer.

Only a few mathematical models have been developed to describe the time-dependent  calcium dynamics within the sperm cell \citep{Li14,olson_10,simons_18,Wennemuth_03}. \cite{Wennemuth_03} studied in detail the calcium clearance phenomena along the flagellum via ATP-ase pumps. Later, \cite{olson_10} developed a model that couples the ATP-ase pumps with the contributions of CatSper channels and a calcium store in the neck, \revtwo{motivated by experimental recordings of calcium fluorescence by \cite{xia_07}}. Of the many modeling studies related to sperm motility, only a handful of models have attempted to \revtwo{consider the role of calcium dynamics on flagellar bending}. Models \revtwo{have incorporated  the effect of calcium indirectly by prescribing} different waveforms for an activated (low calcium, \revtwo{symmetric waveform}) or hyperactivated (high calcium, \revtwo{asymmetric waveform}) sperm\revtwo{, without including any calcium dynamics directly in the model,} to investigate emergent trajectories and interactions with a wall \citep{curtis_12,Ishimoto_16,Ishimoto_17,olson_15}. \revtwo{Another modeling approach has been to directly couple the sperm motility to either a temporal or spatiotemporal evolving calcium concentration via an asymmetric} calcium dependent curvature model where flagellar bending was planar (2D) \citep{olson_11,olson_13b,simons_14}. \revtwo{We remark that \cite{olson_11} and \cite{simons_14} explored both a spatiotemporal calcium coupling to a preferred beat form as well as a preferred beat form with a constant asymmetry. In these models, the evolving calcium dynamics were able to better match experimental data for emergent waveforms and trajectories. In the case of \cite{simons_14}, the direct calcium coupling to preferred amplitude also led to swimming trajectories that allowed escape from a planar wall, matching experimental observations of \cite{chang_12}. Thus, accounting for the relevant spatiotemporal evolution of calcium within the sperm flagellum is important to accurately model sperm motility.}

Here, we develop the first mathematical model that couples the 3D dynamics of the sperm flagellum and the surrounding fluid, modeled respectively as a Kirchhoff rod and a Netwonian viscous fluid, with the CatSper channel mediated calcium dynamics inside the flagellum, via a curvature dependent coupling. The model is then used to investigate the emergent 3D waveforms and trajectories when coupling calcium and curvature, comparing to the 2D case. The planar bending case is fully characterized for the Kirchhoff rod model and compared to previous results when using an Euler elastica representation. Further, we investigate helical bending in the case of equal bending amplitudes (spiral bending) and the case of unequal bending amplitudes (quasi-planar since one amplitude is significantly smaller). The quasi-planar and helical bending cases exhibit emergent trajectories that can be described as a hypotrochoid, similar to the  flagelloid curve observed in experiments by \cite{Woolley03} (and reproduced in Figure \ref{exppic}(a)).

%%%%%%%%%%%%%%%%%%%%%%%%%%%%%
\section{Methods}\label{sec:methods}
%%%%%%%%%%%%%%%%%%%%%%%%%%%%%
\subsection{Flagellum and fluid dynamics}\label{KRfluid}%Flagellum dynamic}
\begin{figure}[!tb]
\centering
\includegraphics[width=0.8\linewidth]{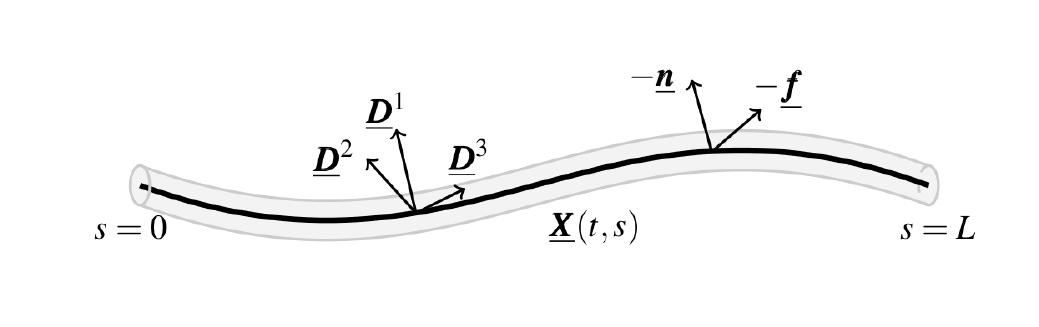}
\caption{\textit{Sketch of the Kirchhoff rod model for the flagellum.} The rod (gray) is represented by the space curve $\vector{X}(t,s)$ (thick black line) and an associated orthonormal triad $\{\vector{D}^1(t,s),\vector{D}^2(t,s),\vector{D}^3(t,s) \}$, $s$ is the rod spatial coordinate, $L$ is the length of the unstressed rod, and $t$ is the temporal coordinate. $\vector{f}(t,s)$ and $\vector{n}(t,s)$ are the external force and torque per unit of length applied to the rod.}
\label{fig:rod}
\vspace*{-9pt}
\end{figure}

Since we are interested in studying the motion of sperm in 3D, we utilize a Kirchhoff rod with preferred curvature and twist to model the elastic flagellum, as in \cite{olson_13}. 
As depicted in Figure~\ref{fig:rod}, the flagellum is represented by a 3D space curve $\vector{X}(t,s)$ and an associated orthonormal triad $\{\vector{D}^1(t,s),\vector{D}^2(t,s),\vector{D}^3(t,s) \}$ for $0 \leq s \leq L$, where $L$ is the length of the unstressed rod, $s$ is a Lagrangian parameter initialized as the arc length and $t$ is time. The rod is assumed to have a circular cross-section of constant radius that is much smaller than the rod length $L$.
%%%

The following elastic energy penalty is considered 
\begin{equation}
E(t) = \Frac{1}{2} \dint_0^L \sum_{i = 1}^3 \left( a_i \left(\Frac{\partial \vector{D}^j}{\partial s} \cdot \vector{D}^k -\Omega_i\right)^2 + b_i \left(\Frac{\partial \vector{X}}{\partial s} \cdot \vector{D}^i - \delta_{3i} \right)^2\right)ds,
\label{eq:energy}
\end{equation}
where $(i,j,k)$ is any cyclic permutation of (1,2,3) and $\delta_{ij}$ is the Kronecker delta. The material parameters for the flagellum (Kirchhoff rod) include the bending moduli $a_1$ and $a_2$,  the twisting modulus $a_3$, shear moduli $b_1$ and $b_2$, and extensional modulus $b_3$. We consider an axially symmetric and isotropic rod, hence $a_1=a_2$ \revtwo{and the material parameters are assumed to be constant along the flagellum length.}
The \textit{preferred strain twist vector} is defined as $(\Omega_1,\Omega_2,\Omega_3)$ where $\Omega_1$ and $\Omega_2$ are the geodesic and normal curvature, respectively, 
\begin{equation}\label{eq:curva}
\Omega = \sqrt{\Omega_1^2 + \Omega_2^2}
\end{equation}
is the \textit{preferred curvature}, and $\Omega_3$ is the \textit{preferred twist}. In this framework, the rod will tend to minimize its energy and differences between the rod configuration and its preferred shape will generate force and torque along the centerline. As described in Section \ref{curvamodel}, the preferred strain twist vector will be time and spatially dependent and coupled to the evolving calcium concentration inside the sperm flagellum. In order to study fully 3D movement and couple mechanics to chemical concentrations, we simplify the sperm representation by neglecting the head or cell body, as in previous modeling studies \citep{olson_13,simons_14}. \rev{This reduces computational costs while still obtaining swimming speeds of the right order of magnitude \citep{Ishimoto_15}.} 

Starting from the energy penalty in \eqref{eq:energy}, and using a similar variational argument to the one detailed in \cite{peskin_02} and \cite{lim_08}, the following balance of linear and angular momentum equations can be derived 
\begin{align}
\vector{0} &= \vector{f} + \deriv{}{\vector{F}}{s},\label{eq:fbal}\\
\vector{0} &= \vector{n} + \deriv{}{\vector{N}}{s} + \left( \deriv{}{\vector{X}}{s} \times \vector{F}\right),\label{eq:tbal}
\end{align}
where $\vector{f}$ and $\vector{n}$, illustrated in Figure~\ref{fig:rod}, are the external force and torque per unit of length applied to the rod, respectively, and $\vector{F}$ and $\vector{N}$ are the average internal force and momentum transmitted across a cross-section of the rod, respectively. The constitutive equations for internal force and torque can be expressed as
\begin{align}
\vector{F} &= \sum_{i=1}^3 b_i \left( \Frac{\partial \vector{X}}{\partial s} \cdot \vector{D}^i - \delta_{3i}\right) \vector{D}^i, \label{eq:int_force}\\
\vector{N} &=\sum_{i=1}^3 a_i \left( \Frac{\partial \vector{D}^j}{\partial s} \cdot \vector{D}^k -\Omega_i\right) \vector{D}^i. \label{eq:int_torque}
\end{align}
We remark that, by minimizing the energy formulation considered in~\eqref{eq:energy}, we are \textit{weakly} imposing the following constraints
\begin{enumerate}
\item the \textit{actual strain twist vector} 
$(\Omega^*_1,\Omega^*_2,\Omega^*_3)=\left(\Frac{\partial \vector{D}^2}{\partial s} \cdot \vector{D}^3,\Frac{\partial \vector{D}^3}{\partial s} \cdot \vector{D}^1,\Frac{\partial \vector{D}^1}{\partial s} \cdot \vector{D}^2\right)$
is equal to the \textit{preferred strain twist vector} $(\Omega_1,\Omega_2,\Omega_3)$, hence the \textit{actual curvature}
\begin{equation}\label{eq:act_curva}
\Omega^*=\sqrt{\left(\Omega^*_1\right)^2 + \left(\Omega^*_2\right)^2}=\sqrt{\left(\Frac{\partial \vector{D}^2}{\partial s} \cdot \vector{D}^3\right)^2 + \left(\Frac{\partial \vector{D}^3}{\partial{s}}\cdot \vector{D}^1\right)^2}
\end{equation}
is equal to the \textit{preferred curvature} $\Omega$ in~\eqref{eq:curva};
\item $\vector{D}^3$ is aligned with the tangent vector, both $\vector{D}^1$ and $\vector{D}^2$ are orthogonal to the tangent vector, and the rod is inextensible, i.e $\left|\deriv{}{\vector{X}}{s}\right| =1$.
\end{enumerate}
These constraints are imposed weakly, hence they tend to be maintained approximately instead of exactly. Moreover, $a_i$ and $b_i$, \revtwo{which} physically represent the material properties, \revtwo{such as the bending and shear resistance of the microtubules and outer dense fibers of the sperm flagellum}, act as Lagrange multipliers for these constraints. Note that the first constraint is reflected in the internal torque constitutive equation~\eqref{eq:int_torque}, while the second constraint is reflected in the internal force constitutive equation~\eqref{eq:int_force}. \revtwo{\cite{lim_08} showed that} the standard Kirchhoff rod model, i.e. the strongly constrained model, can be obtained from the weakly constrained model by setting $b_1=b_2=b_3$, \revtwo{therefore the assumption of homogeneity of the shear and extensional moduli $b_i$ is considered in this work.}

Sperm motility occurs in a regime where viscous forces dominate and acceleration is negligible. Hence, the fluid surrounding the flagellum is modeled as a viscous and incompressible Newtonian  fluid using Stokes equations
\begin{align}
\vector{0}&=-\grad p + \mu \nabu^2{\vector{v}} + \vector{f}^r,\label{eq:stokes}\\
0&=\Div{\vector{v}}\label{eq:incomp},
\end{align}
where $p$ and $\vector{v}$ are the fluid pressure and velocity, respectively, $\mu$ is the fluid viscosity and $\vector{f}^r$ is the force per unit of volume that the sperm exerts on the fluid. 

To solve this fluid-structure interaction problem we use the method of regularized Stokeslets, described in detail in \cite{cortez_01} and \cite{olson_13}. The main idea is to derive the fundamental solution of the Stokes problem in the case of a regularized point force or a regularized point torque. Then, the global solution is obtained by adding the various contributions 
along the flagellum, taking advantage of the linearity of the Stokes equations. The point force or torque applied at the point $\vector{X}_0$ is regularized using the following radially symmetric blob function
\begin{equation}
\phi_{\varepsilon}(\vector{x},\vector{X}_0) = \Frac{15\varepsilon^4}{8\pi \left(\left|\vector{x}-\vector{X}_0\right|^2+\varepsilon^2\right)^{7/2}},\label{eq:blob}
\end{equation}
where $\vector{x}$ is any point in the fluid. The blob function approaches the Dirac delta distribution as $\varepsilon \to 0$ and satisfies $\dint_{\mathbb{R}^3} \phi_{\varepsilon}(\vector{x},\vector{X}_0) d\vector{x} =1$. Note that the regularization parameter $\varepsilon$\revtwo{, assumed to be constant along the flagellum length,} determines 
\revtwo{the region where the majority of the force or torque is spread to the fluid \citep{cortez_05}}
and its value can be chosen so that \revtwo{this region}
corresponds to the physical radius of the cross-section of the rod\revtwo{, for additional details see Section~\ref{sec:numerics} and Table~\ref{tab:param}.}

The dynamic coupling between the surrounding fluid flow and the elastic flagellum is expressed by the force $\vector{f}^r$ in~\eqref{eq:stokes}, which by Newton's third law, depends on the rod external force $\vector{f}$ and torque $\vector{n}$. Using the method of regularized Stokeslets, we can write the force $\vector{f}^r$ exerted by the rod on the point in the fluid $\vector{x}$ as 
\begin{equation}\label{eq:fluid_force}
\vector{f}^r(t,\vector{x})=\int_{\Gamma(t)}\left(-\vector{f}(t,s)+\frac{1}{2}\nabla\times(-\vector{n}(t,s))\right)\phi_{\varepsilon}(\vector{x},\vector{X}(t,s))ds,
\end{equation}
where the curve $\Gamma(t)= \vector{X}(t,s)$ (for more details see \cite{olson_13}).
Note that in~\eqref{eq:stokes}-\eqref{eq:incomp} we consider the incompressible steady Stokes equation, however the \revtwo{force} $\vector{f}^r$ that the rod exerts on the fluid in~\eqref{eq:fluid_force} is time-dependent. 
The kinematic coupling between the surrounding fluid flow and the motion of the elastic flagellum is imposed by the following no-slip boundary conditions on the fluid linear velocity $\vector{v}$ and angular velocity $\vector{w}$
\begin{equation}\label{eq:noslip}
\frac{\partial \vector{X}}{\partial t}=\vector{v}(\vector{X}),\hspace{0.5cm}\frac{\partial \vector{D}^i}{\partial t}=\vector{w}\times\vector{D}^i,\hspace{.2cm}i=1,2,3.
\end{equation}
Accordingly, at each instant in time, given the force exerted by the rod on the fluid, we can solve for the resulting fluid flow and update the rod location assuming it moves with the local fluid velocity.
%%%%%%%%%%%%%%%%%%%%%%%%%%%%%
\subsection{Calcium and curvature dynamics}\label{curvamodel}
In mammalian sperm, the asymmetry and magnitude of flagellar bending is known to vary as a function of the calcium concentration inside the flagellum \citep{HoSuarez01a,Ho02,Lindemann88,Marquez07,Tash82}. We will use a previously developed one-dimensional reaction-diffusion model to account for the relevant calcium dynamics \citep{olson_10}.
The increase in intracellular calcium concentration is initiated by the opening of CatSper channels \revtwo{\citep{Carlson03,xia_07}}, which allows calcium to flow from the external fluid bath to the inside of the sperm flagellum. The sperm is composed of five pieces: the head, neck, mid-piece, principal piece, and end piece~\citep{cummins_85,Pesch_06}. Channels such as CatSper \citep{Chung14} and pumps such as the calcium ATP-ase pump \citep{Okunade} are localized along the length of the principal piece, whereas the Redundant Nuclear Envelope (RNE), a calcium store \citep{HoSuarez01b,Suarez03}, is found in the neck.

The calcium concentration $Ca(t,s)$ at time $t$ and location $s$ along the Kirchhoff rod $\vector{X}(t,s)$ is governed by \begin{equation}
\begin{array}{rcl}
\Frac{\partial}{\partial t} \left( Ca(t,s) \left|\Frac{\partial \vector{X}(t,s)}{\partial s}\right|\right)
&=&D_{Ca} \Frac{\partial}{\partial s}\left(\Frac{\partial Ca(t,s)}{\partial s} \Big{/} \left|\Frac{\partial \vector{X}(t,s)}{\partial s}\right| \right) + J(t,s,Ca(t,s)) \left|\Frac{\partial \vector{X}(t,s)}{\partial s}\right|, 
\label{eq:ca}
\end{array}
\end{equation}
where $D_{Ca}$ is the calcium diffusion coefficient. 
\revtwo{Since we simplify the sperm representation by neglecting the head or cell body, in our model we divide the flagellum in the following pieces: proximal (very small region at the base of the tail in the direction of swimming), neck, mid-piece, principal piece, and end piece. The ranges for these regions are reported in Table~\ref{tab:param}.}
The calcium flux $J$ varies along the sperm length and incorporates the fluxes in the principal piece (CatSper channels and pumps) and in the neck (RNE contribution), depends non-linearly on $Ca$, and is coupled to the evolving concentration of inositol 1,4,5-trisphosphate or $IP_3$  (for more details see \cite{olson_10,olson_11}). Note that \revtwo{$J=0$ in the proximal region, mid-piece and end piece.} Additionally, we assume a large and constant calcium concentration in the fluid surrounding the sperm. The motion and deformation of the rod is taken into account in \eqref{eq:ca} via the term $\left|\partial \vector{X}(t,s)/{\partial s}\right|$,  representing the Jacobian of the transformation from the straight rod to the current rod configuration (for more details see \cite{lai_08} and \cite{stone_90}). 

\begin{figure}[!t] 
\centering
\includegraphics[width=0.9\linewidth]{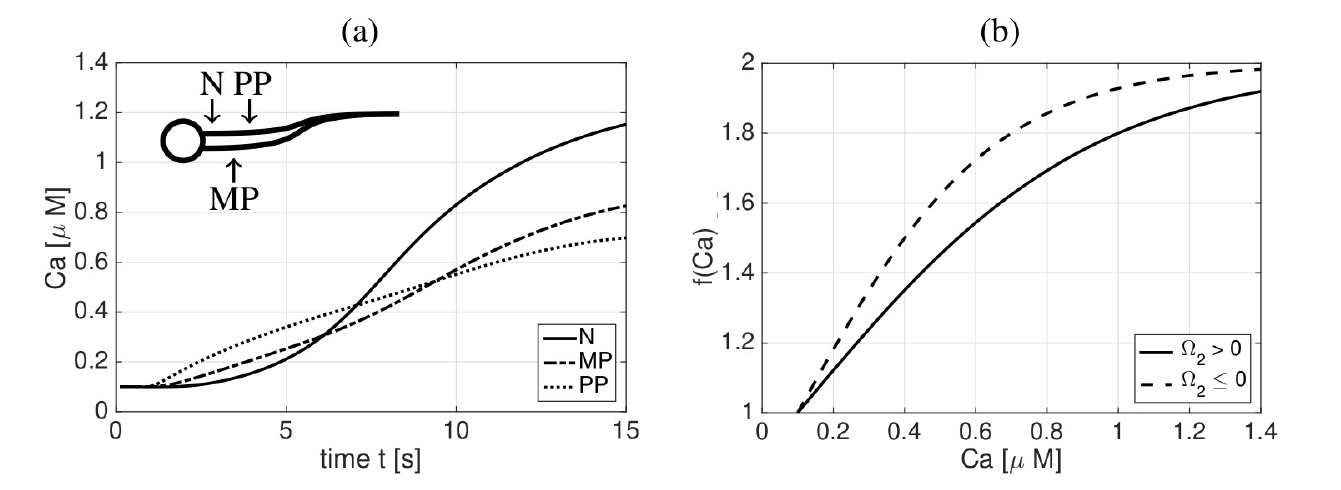}
\caption{\textit{Calcium dynamics.} (a) Calcium ($Ca$) concentration  dynamics in time for the case of a fixed, straight rod at three different points along the sperm length $L$: \revtwo{neck (N) -  at $1\%$ of $L$}, in the mid-piece (MP) - at $16\%$ of $L$, and in the principal piece (PP) - at $24\%$ of $L$. (b) Plot of $f$, given in \eqref{eq:funcf}, as a function of calcium concentration $Ca$ depending on the sign of the preferred normal curvature $\Omega_2$.}
\label{fig:ca}
\vspace*{-9pt}
\end{figure}

Figure~\ref{fig:ca}(a) shows the calcium dynamics in time for the case of a fixed straight rod (i.e. $\left|\partial \vector{X}(t,s)/{\partial s}\right|=1$) at three different points along the the sperm length $L$: \revtwo{neck (N) -  at $1\%$ of $L$}, in the mid-piece (MP) - at $16\%$ of $L$, and in the principal piece (PP) - at $24\%$ of $L$.
\revtwo{The locations of the points N, MP and PP have been chosen to illustrate calcium entry through CatSper channels in the principal piece at earlier time points, calcium induced calcium release in the neck at later time points, and an increase in calcium in the midpiece and principal piece due to diffusion along the length of the flagellum. Note that these calcium trends match the calcium fluorescence results of \cite{xia_07}. Additionally,  we have utilized the model in \cite{olson_10} where we have only updated the scaling of the flagellum to be characteristic of human sperm, as reported in Table~\ref{tab:param}, and we are now neglecting the head region in this work (in contrast to \cite{olson_10}).}
%%%%%

\revtwo{Since the spatiotemporal dynamics will lead to emergent swimming speeds and trajectories, we will couple the movement of the flagellum  to the evolving calcium dynamics through the preferred amplitude, in the same fashion as \cite{olson_11} and \cite{simons_14}. However, in this framework,} we allow for a fully 3D preferred beat form with the Kirchhoff rod model, whereas the previous work used an Euler elastica model with strictly planar or quasi-planar beating. We assume the following preferred reference configuration of the flagellum in time
\begin{equation}
\widehat{\vector{X}}(t,s)=\left[x, \; y, \; z\right]= \left[
x(t,s),  \; A \sin(k \ x(t,s)-\sigma t), \; B \cos(k \ x(t,s)-\sigma t) \right]^{T},\label{eq:wave}
\end{equation}
with the following preferred wave parameters: $A$ and $B$ are the wave amplitudes, $k$ is equal to $2\pi/$ the wavelength, and $\sigma/2\pi$ is the beat frequency. 
Note that depending on the values of $A$ and $B$, the wave can be planar, quasi-planar, or helical \rev{(shown in Figure \ref{fig:3dsperm}).} This is representative of the range of beat forms observed in experiments \citep{Dresdner81,Drobnis88,Guerrero11,Ooi14,Suarez92,Woolley01}. 
\rev{Similar to previous models \citep{Ishimoto_16,Omori16,Smith_09b}, we assume that the two wave amplitudes are related via the chirality parameter $\alpha$ where $B=\alpha A$ for $0\leq \alpha \leq 1$. Here, $\alpha=0$ corresponds to the planar case and $\alpha=1$ corresponds to the helical case.} 
We remark that this assumed reference configuration captures different beat forms but does not try to exactly model flagellar bending at the level of dynein activation or microtubule sliding \citep{Vernon04,Woolley10}. 
\revtwo{Moreover, the preferred reference configuration chosen in~\eqref{eq:wave} causes the rod trajectories to be left-handed, i.e. trajectories show a counterclockwise screwing motion in the direction of swimming.
This is in agreement with experiments, since different species of sperm, such as human and bull sperm, present mostly left-handed trajectories as reported in \cite{Ishijima92}.}

Since increases in calcium are associated with larger amplitude bending \citep{HoSuarez01a,Lindemann88,Marquez07,Tash82}, we assume that the rod amplitudes $A$ and $B$ depend on the calcium concentration along the rod as follows
\begin{equation}
A(t,s)=A_0f(Ca(t,s)),\hspace{0.3cm}\mbox{ and/or}\hspace{0.3cm} B(t,s)=B_0f(Ca(t,s)),\label{eq:amp}
\end{equation}
where $A_0$ and $B_0$ are the constant baseline amplitudes and 
\begin{equation}\label{eq:funcf}
f(Ca(t,s)) = \Frac{2}{1+\exp\left(-\Frac{c_1(Ca(t,s)-\widehat{Ca})}{c_2-\widehat{Ca}} \right)}, 
\end{equation}
with $\widehat{Ca}$ representing the baseline calcium concentration in the flagellum. The function $f(Ca)$ is equal to $1$ at the baseline, i.e. when $Ca=\widehat{Ca}$, and asymptotically approaches  $2$ as $Ca$ increases, see Figure~\ref{fig:ca}(b). Hence, the values of $A$ and $B$ can at most double their baseline values $A_0$ and $B_0$, respectively, depending on the calcium concentration, in agreement with experiments \citep{Marquez07b}.

In terms of the Kirchhoff model presented in Section \ref{KRfluid}, we need to define a preferred strain twist vector that will define flagellar bending and is also coupled to the calcium concentration. Using the reference or preferred beat form $\widehat{\vector{X}}$ in \eqref{eq:wave}, we specify the reference orthonormal triad as $\widehat{\vector{D}}^3$ as the tangent, $\widehat{\vector{D}}^1$ as the normal, and $\widehat{\vector{D}}^2$ as the binormal with respect to $\widehat{\vector{X}}$ \rev{(details given in Appendix \ref{AppCurva})}. Then, 
we can calculate the \textit{preferred strain twist vector} $(\Omega_1,\Omega_2,\Omega_3)$ using 
$\Omega_i = {\partial_s \widehat{\vector{D}}^j} \cdot \widehat{\vector{D}}^k$, 
where $(i,j,k)$ is a cyclic permutation of $(1,2,3)$. 
For this rod configuration, the corresponding strain twist vector is given as  
\begin{subequations}\label{equ:Omega1-3}
\begin{align}
\Omega_1 \revtwo{(t,s)}&= \frac{ B k^2 (1+A^2k^2)\cos ( k \ x(t,s) - \sigma t)} {\sqrt{\mathcal{K}_c} (\mathcal{K}_c+ \mathcal{K}_s  )^{3/2} },\\%\hspace{0.3cm}
\Omega_2 \revtwo{(t,s)}&= -\frac{ A k^2 \sin ( k \ x(t,s) - \sigma t)} {\sqrt{\mathcal{K}_c} ( \mathcal{K}_c+  \mathcal{K}_s) }, \\%\hspace{0.3cm}
\Omega_3 \revtwo{(t,s)}&= \frac{ A B k^3 \sin^2 ( k \ x(t,s) - \sigma t)} {(\mathcal{K}_c) ( \mathcal{K}_c+  \mathcal{K}_s ) },
\end{align}
\end{subequations}
where
\begin{equation}
\mathcal{K}_c=1 + A^2 k^2  \cos^2( k \ x(t,s) - \sigma t)\hspace{0.25cm}\mbox{and}\hspace{0.25cm}\mathcal{K}_s=B^2 k^2 \sin^2(k \ x(t,s) -\sigma t).
\end{equation}
The rod is initialized from \eqref{eq:wave} at $t = 0$ \rev{(additional details given in Appendix \ref{AppCurva})} and it is given a \revtwo{spatiotemporal} evolving strain and twist vector by~\eqref{equ:Omega1-3}; the amplitude in \eqref{eq:amp} also varies \revtwo{in time and space,} and is coupled to the evolving calcium concentration. 
\revtwo{As shown in previous studies \citep{Lim10}, if a helix was given constant preferred curvature and twist (with $\Omega_3\neq0$ and at least one of $\Omega_i\neq0$ for $i=1,2$), it would stop moving once it reaches its preferred helical configuration. However, in this model, since preferred curvature and twist ($\Omega_i$ for $i=1,2,3$) are functions of both space and time, as given in~\eqref{equ:Omega1-3}, the rod will continue to move and propagate a helical wave. We remark that the motion of a sperm flagellum, as opposed to bacterial flagella motion, is not driven by the rotation of a motor at its base~\citep{Park17}, but is driven by dyneins and local force generation along the entire flagellum. The local action of the dyneins is represented in the model as a} propagating curvature wave, similar to observations of human sperm flagellar bending \citep{smith_09}. 

Additionally, as calcium concentrations increase, experiments have shown that flagellar bending shows an asymmetry: the bend grows in the principal bend direction and ends early in the reverse bend \citep{Tash82}. We include this asymmetry in our model by having the parameter $c_2$ in~\eqref{eq:funcf} vary the speed at which the function $f$ approaches its maximum, as shown in Figure~\ref{fig:ca}(b), coupling $c_2$ to the sign of the curvature ${\Omega}_2$ as follows
\begin{equation}\label{eq:c2}
c_2=\left\{
\begin{array}{cc}
c_{2,p}& \mbox{if }{\Omega}_2(t,s)>0, \\ 
c_{2,n}& \mbox{if }{\Omega}_2(t,s)\leq0,
\end{array}
\right. 
\end{equation}
similar to \cite{olson_11} and \cite{simons_14}. In order to introduce an asymmetry in flagellar bending, both in the case of a 3D wave (i.e. $A\neq0$ and $B\neq0$), and in the case of a planar wave (i.e. $B=0$), we vary the parameter $c_2$ in~\eqref{eq:c2} based on the sign of ${\Omega}_2$, which is the only component of the reference preferred strain vector in~\eqref{equ:Omega1-3} that does not depend on $B$.

\subsection{Numerical algorithm for coupling}\label{sec:numerics}
At time $t=0$, the rod and orthonormal triads are initialized using the preferred reference configuration in \eqref{eq:wave} for a given set of wave parameters and \revtwo{the rod} is discretized into \revtwo{$M+1$} points with constant spacing $\triangle s=\revtwo{L/M}$\revtwo{, i.e. $s_k=k \triangle s$ for $k=\revtwo{0},\ldots,M$} \rev{(additional details given in Appendix \ref{AppCurva})}. The following steps are followed to solve for the new configuration of the rod at time $\triangle t$.
\begin{enumerate}
\item  Compute the calcium concentration along the rod by solving \eqref{eq:ca} as in \cite{olson_10,olson_11}. A symmetric Crank-Nicolson scheme detailed in \cite{lai_08} is used, which ensures that the total mass of calcium is numerically conserved in the case of flux $J=0$;
\item Given the calcium concentration, determine the preferred amplitude in \eqref{eq:amp} for $A$ and $B$ at each discretized point along the rod;
\item Determine the preferred strain and twist $\Omega_i$ in \eqref{equ:Omega1-3} for each of the \revtwo{$M+1$} points, assuming $A$ and $B$ constant at each point;
\item Calculate the point forces $\vector{f}_k$ and torques $\vector{n}_k$ along the rod using \eqref{eq:fbal}-\eqref{eq:int_torque} as in \cite{olson_13};
\item Solve for the fluid flow in \eqref{eq:stokes}-\eqref{eq:incomp} using regularized fundamental solutions: the local linear velocity $\vector{v}$ and angular velocity $\vector{w}$ are given as
\begin{align}
\vector{v}(\vector{x})&=\frac{1}{\mu}\sum_{\revtwo{k=0}}^M\mathcal{S}[-\vector{f}_k\triangle s]+\frac{1}{\mu}\sum_{\revtwo{k=0}}^M\mathcal{R}[-\vector{n}_k\triangle s],\\
\vector{w}(\vector{x})&=\frac{1}{2}\nabla\times\vector{v}=\frac{1}{\mu}\sum_{\revtwo{k=0}}^M\mathcal{R}(-\vector{f}_k\triangle s)+\frac{1}{\mu}\sum_{\revtwo{k=0}}^M\mathcal{D}(-\vector{n}_k\triangle s), \end{align}
where $\mathcal{R}$, $\mathcal{S}$, and $\mathcal{D}$ are the regularized rotlet, Stokeslet, and dipole for the blob function given in \eqref{eq:blob}, as detailed in \cite{olson_13};
\item Update the rod location $\vector{X}$ and orthonormal triads $\vector{D}^i$ for $i=1,2,3$ using the local linear and angular velocity via the no-slip conditions in~\eqref{eq:noslip}. The forward Euler method 
is used to solve these equations.
\end{enumerate}
%%%%%%%%%%%%%%%%%%%%%%%%%%%%%
\section{Results and Discussion}\label{sec:results}
\begin{figure}[!t]
\centering
\includegraphics[width=0.7\linewidth]{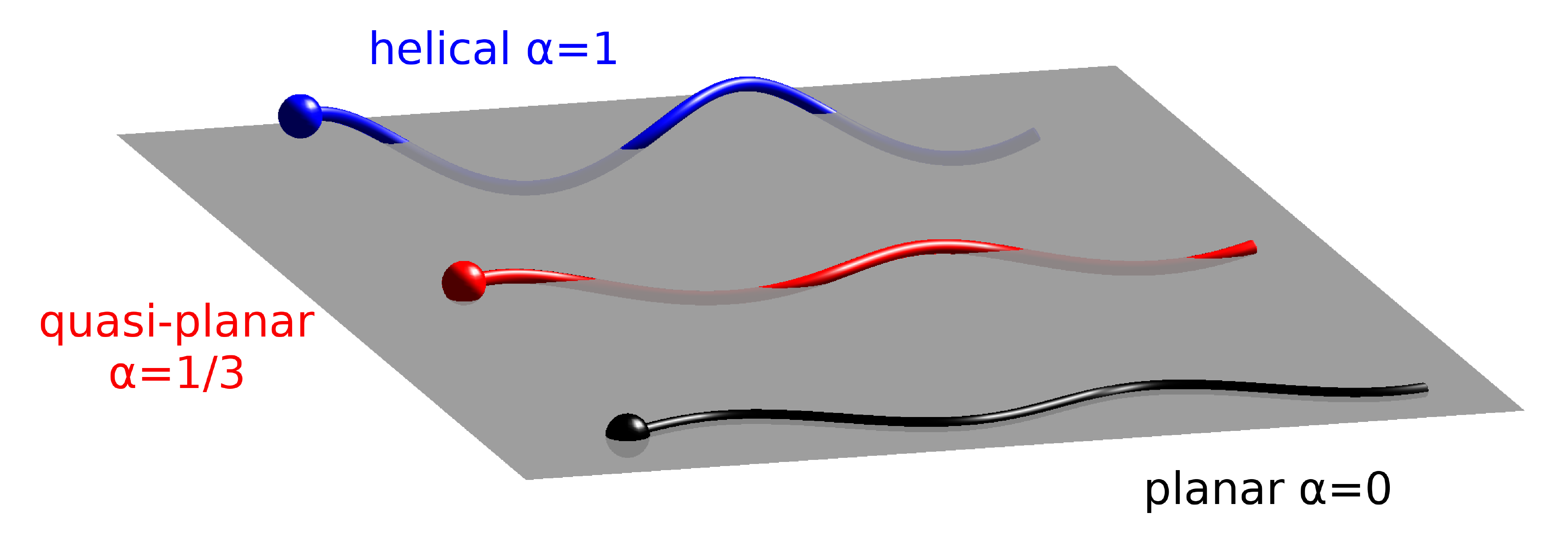}
\caption{\rev{\textit{Flagellar beat forms.} Schematic representation of the three flagellar beat forms considered: planar wave (bottom, black), quasi-planar wave (middle, red) and helical wave (top, blue), and their corresponding value of the chirality parameter $\alpha$.}}
\label{fig:3dsperm}
\vspace*{-9pt}
\end{figure}

We investigate the effect of the calcium and curvature coupling described in Section~\ref{curvamodel} on sperm motility. The flagellar beat forms we consider are illustrated in Figure~\ref{fig:3dsperm} and are given as follows,
\begin{itemize}
\item 2D dynamics: \textit{planar wave} with $A_0=3$ $\mu$m and $B_0=0$ $\mu$m , i.e. $\alpha=0$;
\item 3D dynamics: 
\begin{itemize}
\item[\textit{i)}] \textit{helical wave} with $A_0=B_0=3$ $\mu$m, i.e. $\alpha=1$;
\item[\textit{ii)}] \textit{quasi-planar wave} with $A_0=3$ $\mu$m and $B_0=1$ $\mu$m, i.e. $\alpha=1/3$.
\end{itemize}
\end{itemize}
Moreover, we consider four possible cases of calcium and curvature coupling
\begin{enumerate}
\item[\textit{a)}] no-coupling (\textit{No Ca}): $A$ and $B$ are fixed to the baseline values, i.e. $A=A_0$ and $B=B_0$;
\item[\textit{b)}] symmetric coupling (\textit{Ca sym}): $A$ and $B$ vary symmetrically with respect to curvature, i.e. in ~\eqref{eq:funcf} $c_2=c_{2,p}=c_{2,n}=1\mu$M; 
\item[\textit{c)}] asymmetric coupling (\textit{Ca asym $A$ \& $B$}): $A$ and $B$ vary asymmetrically with respect to curvature, i.e. in ~\eqref{eq:funcf} $c_{2,p}\neq c_{2,n}$ as in Table~\ref{tab:param};
\item[\textit{d)}] asymmetric coupling only A (\textit{Ca asym $A$}): $A$ varies asymmetrically with respect to curvature, while $B$ is kept constant to its baseline value, i.e. $B=B_0$.
\vspace{-4pt}
\end{enumerate}
Note that in the planar wave case, since $B_0=0$, couplings \textit{c)} and \textit{d)} are equivalent, and we will refer to them as the asymmetric coupling case (\textit{Ca asym}).
%%%%%

The material and geometric parameters of the sperm flagellum\revtwo{, as well as the calcium model parameters,} are summarized in Table~\ref{tab:param}, with references given where applicable.  \rev{In particular, the flagellar beat form parameters $L$, $A_0$, $B_0$, $\sigma$ and $k$, reported in Table~\ref{tab:param}, have been taken directly from the reported references. 
\revtwo{The values for the different regions of the flagellum are reported in Table~\ref{tab:param}. These ranges are based on the average dimensions of the mammalian sperm flagellum reported in~\cite{cummins_85} and \cite{Pesch_06}, following the same proportions between region size and flagellum size used in~\cite{olson_10} and accounting for the fact that in this model, in contrast to~\cite{olson_10}, we are neglecting a sperm head region.}
The values for the material properties of a mammalian sperm flagellum, i.e. $a_i$ and $b_i$, have been determined by starting from the shear and bending resistance values extracted from sea urchin sperm in~\cite{Pelle09} and \cite{okuno_79}. These were adjusted by one or two orders of magnitude to take into account the fact that mammalian sperm are generally stiffer than invertebrate sperm due to the presence of the outer dense fibers~\citep{Lindemann73,Schmitz04,Zhao18}. The values of $D_{Ca}$ and $\widehat{Ca}$ have been taken directly from the reported references, in accordance to~\cite{olson_10}, while the values of $c_1$, $c_{2,p}$ and $c_{2,n}$ have been chosen phenomenologically. For example, $c_1$ is chosen so that the function $f(Ca)$, depicted in Figure~\ref{fig:ca}(b), attains  $90\%$ of its upper bound for $Ca=c_2$. Considering that the calcium dynamics, reported in Figure~\ref{fig:ca}(a), \revtwo{are} fully developed after $t=10$s, we choose  $c_{2,n}=1\mu$M, and $c_{2,p}=0.7\mu$M$<c_{2,n}$ to introduce the bending asymmetry by means of changing the slope of the function $f$.}

\revtwo{The values for the numerical algorithm parameters are also listed in Table~\ref{tab:param} and were chosen to ensure convergence and stability of computational results. In particular, the value of the regularization parameter $\epsilon$ is chosen to comply to its physical and numerical interpretation. We note that the singular Stokeslet solution is obtained in the limit as $\epsilon\to 0$ and the magnitude of $\epsilon$ will in turn control the regularization error. Thus, we want to minimize error but still have a physically accurate representation of the flagellum since, via the blob function~\eqref{eq:blob}, $\epsilon$ also represents the radius of the spherical region around a point where the majority of the force or torque is spread to the fluid. This is why we can also view $\epsilon$ as a physical parameter: the virtual radius of the flagellar cross-section. We chose $\epsilon=1\mu$m in agreement with the characteristic radii for the mid-piece of mammalian sperm reported in~\cite{cummins_85}.
Moreover, in order for the derivation of force and torque balance on the Kirchhoff rod to be valid, we must have that the rod radius is much smaller than the length, and this is satisfied with the reported virtual rod radius $\epsilon$ and rod length $L$. 
The proportionality between $\triangle s$ and $\epsilon$ is chosen to make sure that the blob functions overlap sufficiently and to ensure that the fluid does not leak through the rod centerline~\citep{cortez_18}. The value of $\triangle t$ is chosen small enough to ensure numerical stability of the forward Euler method used for the time discretization~\citep{olson_13}.} %

In this section, we omit results for calcium concentration dynamics along the flagellum since they are analogous to the case of a fixed rod shown in Figure~\ref{fig:ca}(a), with the minor addition of infinitesimal oscillations due to the deformation of the rod, see \cite{olson_11} for more details.
\begin{table}[!t]
{
\scalebox{0.9}{
\mbox{
\tabcolsep=8pt\begin{tabular}{@{}llll@{\hskip -0.2in}l@{}}
%\tblhead{
\toprule
{\bf Parameter }& {\bf Value} & {\bf Units} &  \multicolumn{2}{l}{\bf Reference}\\%[-9.5pt]}\\[-9.5pt]
\toprule
Amplitude range, $A_0$ and $B_0$ & $[0,3]$ & $\mu$m & \revtwo{\multirow{3}{*}{$ \ \left\} \begin{array}{c} \\ \\ \\ \end{array} \right.$}} &  \multirow{3}{*}{ \cite{smith_09}}\\ 
Beat frequency, $\sigma$ & $20(2\pi)$ & Hz & \\
Wavelength, $2\pi/k$ & $30$& $\mu$m & \\
\arrayrulecolor{black!30}\midrule
 Length, $L$ & $60$&$\mu$m & \revtwo{\multirow{6}{*}{ $\left\} \begin{array}{c} \\ \\ \\ \\ \\ \\ \end{array} \right.$} }& \multirow{6}{*}{ \makecell[l]{\cite{cummins_85} \\ \cite{Pesch_06}}}\\
\revtwo{Proximal range} & $[0,0.5]\% L$ & $\mu$m & \\
Neck range & $[0.5,2.5]\% L$ & $\mu$m &\\
Mid-piece range & $[2.5,20.5]\% L$ & $\mu$m &\\
Principal piece range & $[20.5,93]\%L$ & $\mu$m & \\
End piece range & $[93,100]\%L$ & $\mu$m & \\
\midrule
Fluid viscosity, $\mu$ & $1\times10^{-6}$ & g $\mu$m$^{-1}$ s$^{-1}$ &  \multicolumn{2}{l}{(water at room temperature)}\\
\midrule
Bending moduli, $a_1=a_2$ & $1$ & g $\mu$m$^{3}$ s$^{-2}$ & \revtwo{\multirow{4}{*}{ $\left\} \begin{array}{c} \\ \\ \\ \\ \end{array} \right. $} }& \multirow{4}{*}{ \makecell[l]{\cite{Pelle09}\\ \cite{okuno_79}}}\\
Twist modulus, $a_3$ & $1$ & g $\mu$m$^{3}$ s$^{-2}$ & \\
Shear moduli, $b_1=b_2$ & $0.6$ & g $\mu$m s$^{-2}$ &\\
Extension modulus, $b_3$ & $0.6$ & g $\mu$m s$^{-2}$ & \\
\midrule
\multirow{2}{*}{Calcium diffusion coefficient, $D_{Ca}$ }& \multirow{2}{*}{$20$} & \multirow{2}{*}{$\mu$m$^{2}$ s$^{-1}$} &  \multicolumn{2}{l}{\cite{Allbritton92}} \\
 & & &  \multicolumn{2}{l}{\cite{sneyd_95}}\\
 \arrayrulecolor{black!30}\midrule
Calcium baseline concentration, $\widehat{Ca}$ & $0.1$ & $\mu$M & \multicolumn{2}{l}{\cite{Wennemuth_03}}\\
\midrule
Calcium parameter, $c_{1}$& $\ln(9)$& 1&\\
Calcium parameter, $c_{2,p}$& $0.7$& $\mu$M &\\
Calcium parameter, $c_{2,n}$& $1$& $\mu$M&\\
\midrule
\multirow{2}{*}{Rod spatial discretization, $\triangle s$} & \multirow{2}{*}{$0.2$} & \multirow{2}{*}{$\mu$m} &   \multicolumn{2}{l}{\revtwo{\cite{lee_14}}}\\ 
& & & \multicolumn{2}{l}{\revtwo{\cite{olson_13}}}\\
\arrayrulecolor{black!30}\midrule
Temporal discretization, $\triangle t$ & $1\times10^{-6}$ & s &  \multicolumn{2}{l}{\revtwo{\cite{olson_13}}}\\
\arrayrulecolor{black!30}\midrule
Regularization parameter, $\varepsilon$ & \multirow{2}{*}{$5\triangle s = 1$} & \multirow{2}{*}{$\mu$m} &  \multicolumn{2}{l}{\revtwo{\cite{cummins_85}}}\\
\revtwo{(Virtual rod radius)} & & & \multicolumn{2}{l}{\revtwo{\cite{olson_13}}}\\
    \arrayrulecolor{black}\bottomrule
\end{tabular}
\label{tab:param}
}}
\caption{The material and geometric parameters of the sperm flagellum, the calcium model parameters, as well as the the numerical algorithm parameters. }
}
\end{table}
%%%%%%%%%%%%%%%%%%%%%%%%%%%%%%%%%%%%%
\subsection{2D flagellar dynamics}
\begin{figure}[!t]
\centering
\includegraphics[width=0.92\linewidth]{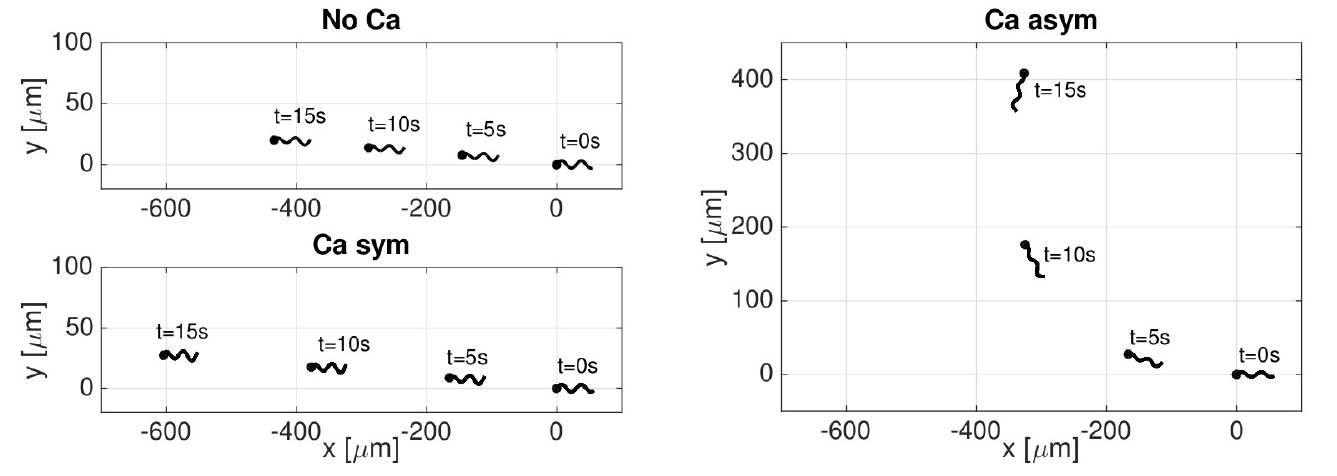}
\caption{\textit{Planar wave.} Sperm trajectories in time $t$ from 0 to $15$s for the case of no calcium-curvature coupling (No Ca), symmetric coupling (Ca sym) and asymmetric coupling (Ca asym). A filled in circle is used to denote the head and swimming direction.}
\label{fig:planar_panel_time}
\vspace*{-9pt}
\end{figure}

Figure~\ref{fig:planar_panel_time} shows sperm trajectories over a 15 second interval for simulations corresponding to the different cases of calcium-curvature coupling. Note that in the planar wave case $\Omega_1=\Omega_3=0$ since $B=0$, hence the preferred curvature of the rod is $\Omega=|\Omega_2|$. For this curvature case, in both the symmetric and asymmetric coupling cases, we observe that the sperm has a higher linear velocity ($\mu$m s$^{-1}$) compared to the no-coupling case: \textit{No Ca} - $28.9$, \textit{Ca symm} - 40.4, and \textit{Ca asym} - $34.9$. The linear velocities extracted from the model results are in agreement with the experimental measurements for human sperm reported in \cite{smith_09}.
However, the symmetric coupling is not enough to reproduce the calcium-dependent turn in the sperm trajectory observed in vivo~\citep{Marquez07b}. Only when the beat asymmetry is included in the coupling, the numerical results are able to reproduce the turning phenomenon. This is similar to previous computational results for a 2D model of a sperm using an Euler elastica representation of the flagellum \citep{olson_11}. For this reason, we will only consider the asymmetric coupling cases in the rest of the results section. \revtwo{In the case of asymmetric calcium coupling, we remark that the trajectory starts as straight but as the calcium increases in the principal piece due to CatSper channels opening, this initiates calcium induced calcium release in the neck region, which then causes an increase along the length of the flagellum due to diffusion (Figure \ref{fig:ca}(a)). In Figure \ref{fig:planar_panel_time}, the increased path curvature due to the turning motion starts at $t=5$s when the calcium is increasing and reaches a stable trajectory for $t=10-15$s when calcium is reaching the maximal values along the flagellum.}
%%%%

\begin{figure}[!t]
\centering
\includegraphics[width=1\linewidth]{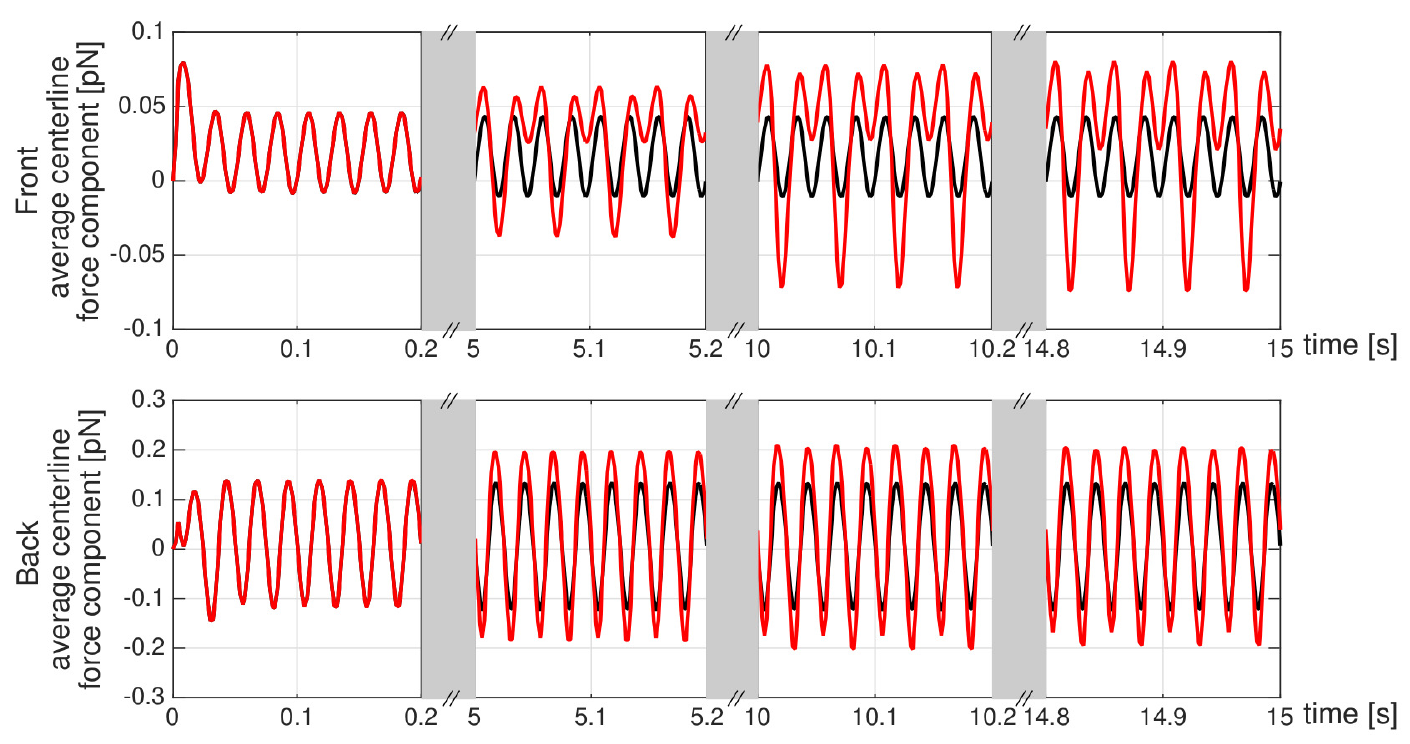}
\caption{\textit{Planar wave.} Average centerline \revtwo{force component $f_{\textit{c}}(t,s)$ in the front region of the rod (first $20$ points of the spatial discretization, top panel) and in the back region of the rod (last $20$ points of the spatial discretization, bottom panel)} for no-coupling (No Ca, black) and asymmetric coupling (Ca asym, red) cases. For $t=0$ to $0.2$s, the curves lie on top of each other.}
\label{fig:planar_force_H_EP}
\vspace*{-9pt}
\end{figure}

To fully understand the dynamics related to the turning phenomenon in the asymmetric coupling case, we compute \revtwo{the centerline component $f_{\textit{c}}(t,s)$ of the} external force applied by the sperm to the fluid. \revtwo{Let $\vector{P}_1(t)$ be the location at time $t$ of the rod first point in the direction of swimming, i.e $\vector{P}_1(t)=\vector{X}(t,s_0)$, and let $\vector{P}_2(t)$ be the location at time $t$ of the rod center of mass, i.e. $\vector{P}_2(t)=\sum_{k=0}^{M} \vector{X}(t,s_k)/(M+1)$. We define the sperm centerline as the line passing through $\vector{P}_1$ and $\vector{P}_2$ and we define the swimming direction as the vector $\vector{u} = \vector{P}_1 - \vector{P}_2$ pointing from the center of mass to the rod first point in the direction of swimming. In this framework, the centerline force component can be expressed as
\begin{equation}
f_{\textit{c}}(t,s) = \Frac{\vector{f}^r (t,s) \cdot \vector{u}(t)}{|\vector{u}(t)|} = \left |\vector{f}^r (t,s) \right | \cos\left({\beta}(t,s)\right),
\end{equation}
where $\vector{f}^r$ is the vectorial force exerted by the rod to the fluid~\eqref{eq:fluid_force} and $\beta$ is the angle between $\vector{f}^r$ and $\vector{u}$. We remark that the sign of the centerline force component $f_{\textit{c}}(t,s)$ depends on the swimming direction $\vector{u}$ via the cosine of the angle $\beta$, namely $f_{\textit{c}}(t,s)$ is negative if the force $\vector{f}^r$ is acting against the swimming direction and is positive otherwise.
Figure~\ref{fig:planar_force_H_EP} shows the time evolution of the average $f_{\textit{c}}(t,s) $ in the front (top panel) and back (bottom panel) regions of the rod for the case of no-coupling and asymmetric coupling. The front region is defined as the first $20$ points of the spatial discretization in the swimming direction, and the back region is defined as the last $20$ points of the spatial discretization.}
\revtwo{After a first transitional region near $t=0$, the average $f_{\textit{c}}$ curves evolve into a symmetric oscillatory behavior in the front and back of the rod for both coupling cases considered.
In the no-coupling case, the symmetric oscillatory behavior does not change with time, giving rise to straight trajectories. In the asymmetric coupling case, the oscillatory behavior with time becomes asymmetric and varies in amplitude and frequency. At time $t=5s$, which corresponds to the turning point (see Figure~\ref{fig:planar_panel_time}), the average $f_{\textit{c}}$ in the front and back regions show an increase in amplitude with respect to the no coupling case, together with a decrease in frequency and the onset of asymmetric oscillations in the front region. In the front, the turning point is followed by a more pronounced asymmetry in the oscillations and a further increase in amplitude for $t=10$s and $t=15$s. In the back, for $t=10$s and $t=15$s, the average $f_{\textit{c}}$ shows a slightly asymmetric oscillatory behavior, while amplitude and frequency do not vary significantly with respect to $t=5$s. This is consistent with the calcium dynamics in the flagellum, reported in Figure~\ref{fig:ca}(a), since the calcium concentration will attain higher values in the front than in the back as time goes by.}
We remark that in the planar case, the sperm trajectories lie in the 2D $xy$-plane, hence the average external torque component along the centerline is zero. In comparison to a 2D model of a sperm using an Euler elastica representation, the forces we present in the direction of the centerline are on the same order of magnitude \citep{simons_14}.
%%%%
%%%
\begin{figure}[!t]
\centering
\includegraphics[width=0.95\linewidth]{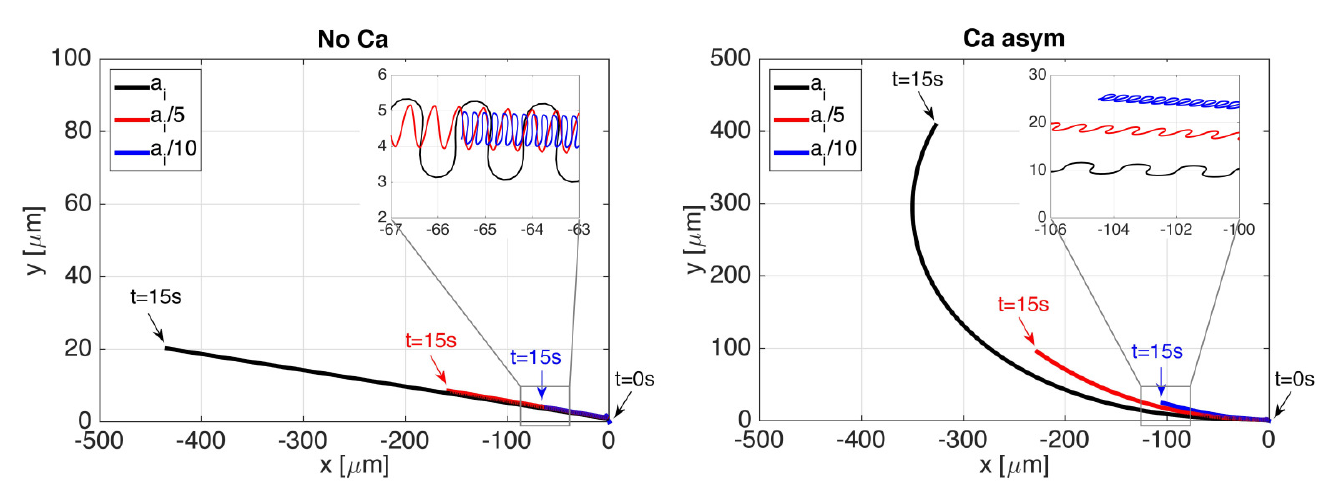}
\caption{\textit{Planar wave.} \revtwo{Trajectories of the flagellum first point in the direction of swimming} in time $t$ varying the bending and twist moduli $a_i$: moduli as in Table~\ref{tab:param} (black), reduced by a factor of $5$ (red) and reduced by a factor of $10$ (blue). No-coupling case (No Ca, left) and asymmetric coupling case (Ca asym, right). The zoomed in portion highlights that the trajectories oscillate in time.}
\label{fig:planar_head}
\vspace*{-9pt}
\end{figure}

\begin{figure}[!t]
\centering
\includegraphics[width=1\linewidth]{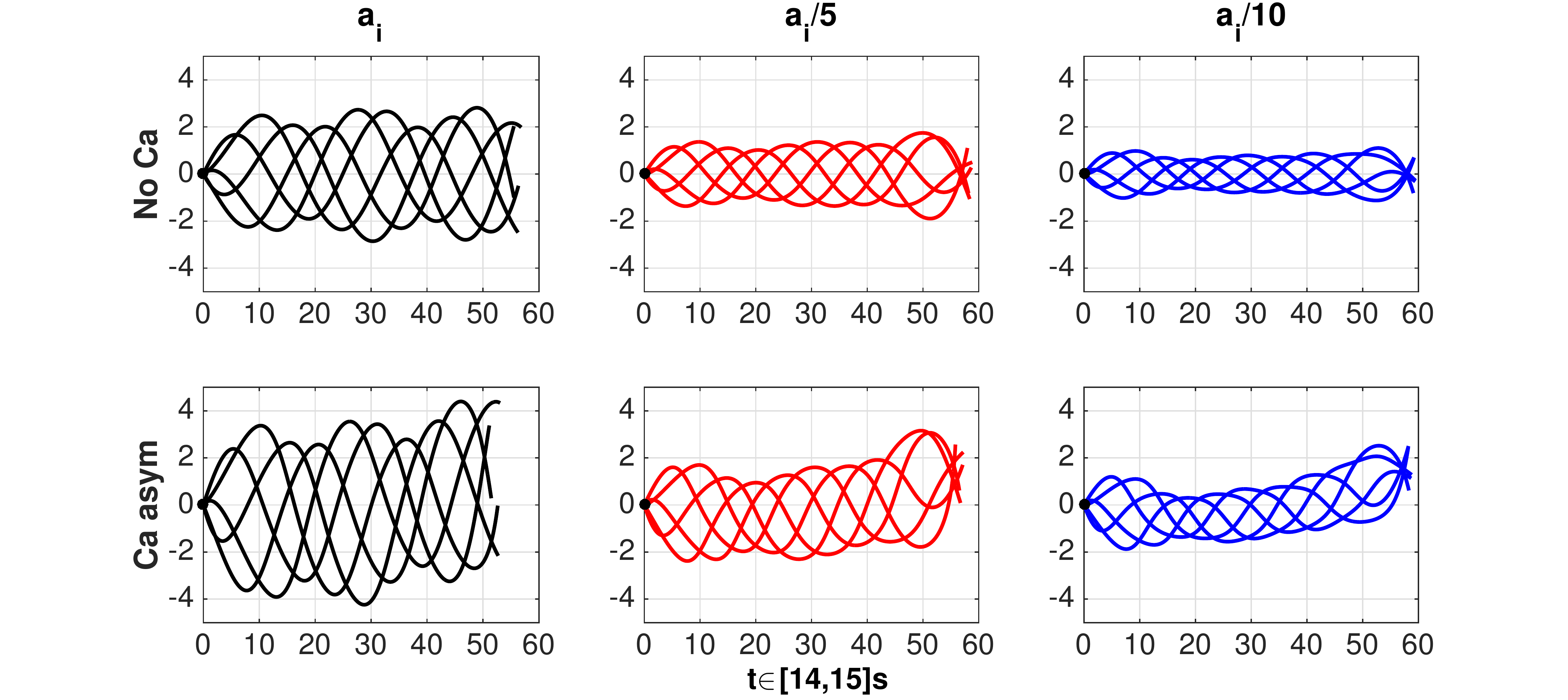}
\caption{\textit{Planar wave.} Evolution of flagellum (translated and rotated to the horizontal axis) for $t=14$s to $15$s, spatial units shown are in microns. No-coupling (No Ca, top row), asymmetric coupling (Ca asym, bottom row), bending and twist moduli as in Table~\ref{tab:param} (left, black), moduli reduced by a factor of $5$ (center, red) and moduli reduced by a factor of $10$ (right, blue). The black circle at the origin denotes the first point on the flagellum.}
\label{fig:planar_tails}
\vspace*{-9pt}
\end{figure}

In Figures~\ref{fig:planar_head} and~\ref{fig:planar_tails} we investigate the effect of varying the flagellum material properties on the motion of the sperm. We can consider the moduli $a_i$ as both numerical and material parameters. As numerical parameters, the moduli are Lagrange multipliers that enforce how strictly the curvature and twist are enforced. In addition, these moduli give an effective stiffness to the elastic rod that represents the sperm flagellum. From experiments, we know that the stiffness varies by several orders of magnitude for different species of sperm, where mammalian sperm are generally stiffer than invertebrate sperm \citep{Lindemann73,Schmitz04,Pelle09}. Thus, we run simulations where the bending and twist moduli $a_i$, in Table~\ref{tab:param}, have been reduced by a factor of $5$ and by a factor of $10$. In this parameter regime, the larger bending and twist moduli correspond to a stiffer elastic rod, \revtwo{i.e. stiffer microtubules and outer dense fibers of the flagellum. 
We remark that in these simulations the radius of the rod cross-section is kept constant, i.e. $\epsilon$ is kept constant.}
Figure~\ref{fig:planar_head} shows the \revtwo{trace of the flagellum first point in the direction of swimming} over 15 seconds in the case of no-coupling (left) and in the case of asymmetric coupling (right) for the three bending and twist moduli values considered. 
\revtwo{For this range of} material properties considered, the \revtwo{rod first point} shows, in general, a linear trajectory in the no-coupling case and a clockwise turning trajectory in the asymmetric coupling case. Moreover, in the asymmetric coupling case, the trajectory radius of curvature is directly proportional to the bending and twist moduli, and the linear speed is also proportional to $a_i$ in both coupling cases considered. While the \revtwo{rod first point} trace shows a general linear or turning trajectory, in reality, as shown in the zoomed area of Figure~\ref{fig:planar_head}, the \revtwo{rod first point} trace is an oscillating curve in time. The amplitude and frequency of these oscillations vary with the material properties and the coupling condition considered, while the shape of the oscillations vary with the material properties but seem to be independent from the coupling condition. 

In Figure~\ref{fig:planar_tails}, the flagellar configurations are illustrated  for different values of bending and twist moduli $a_i$, over the time interval $14$s to $15$s. The configurations are translated and rotated so that the \revtwo{rod first point in the direction of swimming} is at the origin and the centerline lies on the horizontal axis, where spatial units are shown in microns. 
Each row corresponds to a calcium-curvature coupling case, and each column corresponds to different flagellum material properties. The results show that the emergent flagellar wave amplitude, \revtwo{for this range of moduli and preferred beat form parameters, decreases as bending and twisting moduli are decreased. The largest amplitude achieved corresponds to the preferred amplitude. Since the emergent amplitude is partially determined by the flagellum trying to minimize the energy given in \eqref{eq:energy}, the amplitude achieved along the length is non-constant and, as in Figure \ref{fig:planar_tails}, tends to increase from the first point in the direction of swimming to the last point.} For all the values of bending and twist moduli considered, there is an increase in the wave amplitude when the asymmetric coupling is considered, and this increase is proportional to the corresponding no-coupling configuration. Moreover, the asymmetric coupling generates an increased tilt of the end piece with respect to the centerline, and this phenomenon is more prominent as the bending and twisting moduli decrease. We also note that trajectories and observed waveforms for varying moduli in the Kirchhoff rod model are similar to those for an Euler elastica model \citep{olson_11}.
%%%%

%%%%%%%%%%%%%%%%%%%%%%%%%%%%%%%%%%%%%
\subsection{3D dynamics}
%%%%%%%%%%%%%%%%%%%%%%%%%%%%%%%%%%%%%
\subsubsection{Helical wave.}
\begin{figure}[!t]
\centering
\includegraphics[width=0.9\linewidth]{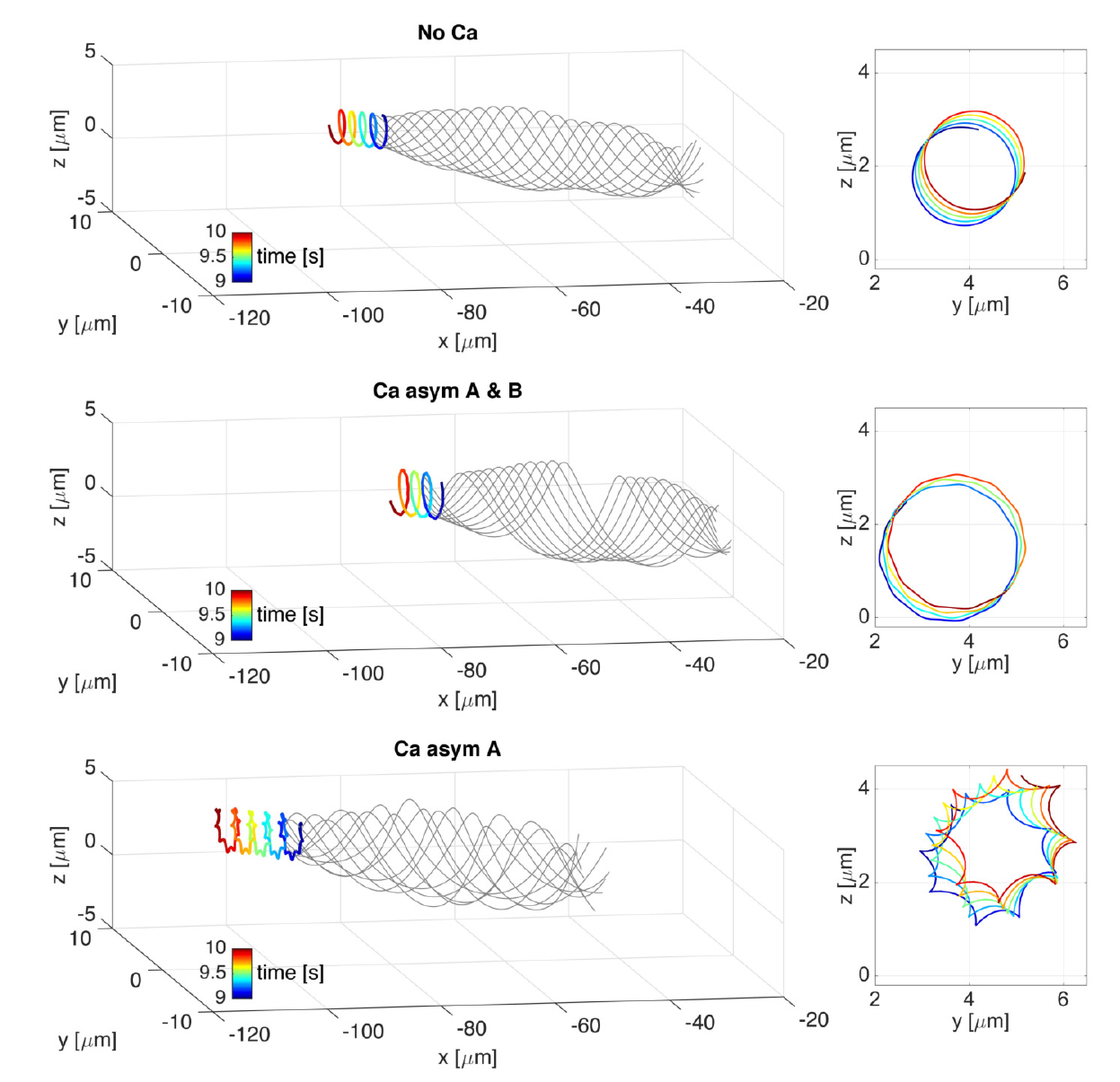}
\caption{\textit{Helical wave.} 3D overlap of flagellar configurations (thin gray lines) and \revtwo{the trajectory of the first point in the direction of swimming} (thick colored line) over the time interval of $1$s, form $9$s (blue) to $10$s (red), on the left. Corresponding flagelloid curve traced by the \revtwo{first point} on the $yz$-plane, on the right. No-coupling case (No Ca, top), asymmetric coupling case (Ca asym $A$ \& $B$, center) and asymmetric coupling only A case (Ca asym $A$, bottom).}
\label{fig:helix_3p}
\vspace*{-9pt}
\end{figure}

In the helical wave case, Figure~\ref{fig:helix_3p} shows the overlap of flagellar configurations (thin gray lines) and \revtwo{the trajectory of the first point in the direction of swimming} (thick colored line) over the time interval of $1$s, form $9$s (blue) to $10$s (red), for three different coupling cases: \textit{No Ca}, \textit{Ca asym $A$ \& $B$} and \textit{Ca asym $A$}. On the right is depicted the corresponding \textit{flagelloid curve} (f-curve) traced by the \revtwo{first point} on the $yz$-plane. 
The f-curve is defined as the path followed by a fixed point on the flagellum~\citep{Woolley98}. In the \textit{No Ca} and \textit{Ca asym $A$ \& $B$} coupling cases, the flagellar configurations remain helical in time, with increased amplitudes in the \textit{Ca asym $A$ \& $B$} case, and the f-curves traced by the \revtwo{first point} are almost circular. In the \textit{Ca asym $A$} coupling case, the flagellar configurations show an irregular beat pattern and the \revtwo{first point} f-curve resembles a hypotrochoid with eight singular points, with increased amplitudes in comparison to the no-coupling case. We note that the linear velocity of the \textit{Ca asym $A$ \& $B$} case is $12\%$ lower than the no coupling case, see Table~\ref{tab:vel}, due to the higher amplitude and almost perfect helical shape, i.e. chirality parameter $\alpha\simeq 1$. While, in the \textit{Ca asym $A$} 
case, the irregularity of the flagellum shape together with the increase in wave amplitudes produces an increase of $53\%$ in the linear velocity compared to the no-coupling case, see Table~\ref{tab:vel}. 
%%%%

\begin{figure}[!t]
\centering
\includegraphics[width=0.9\linewidth]{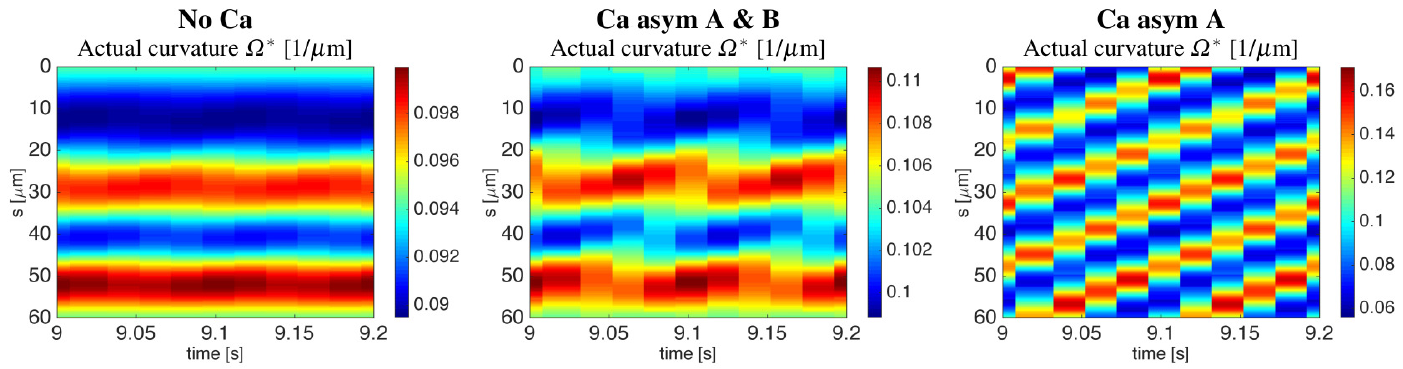}
\caption{\textit{Helical wave.} Actual curvature $\Omega^*$ variations along the flagellum spatial coordinate $s$ in the time interval from $9$s to $9.2$s. No-coupling case (No Ca, left), asymmetric coupling case (Ca asym $A$ \& $B$, center) and asymmetric coupling only A case (Ca asym $A$, right).}
\label{fig:helix_curva}
\vspace*{-9pt}
\end{figure}

Variations in the actual curvature $\Omega^*$~\eqref{eq:act_curva}  along the flagellum spatial coordinate $s$ in the time interval from $9$s to $9.2$s are reported in Figure~\ref{fig:helix_curva}, for the same three coupling cases investigated in~Figure~\ref{fig:helix_3p}. 
In the \textit{No Ca} and \textit{Ca asym $A$ \& $B$} coupling cases the flagellum curvature is almost constant, and this is consistent with the helical flagellum configuration reported in~Figure~\ref{fig:helix_3p}, since helices by definition have constant curvature and twist.
However, in the \textit{Ca asym $A$} case the curvature varies periodically along the flagellum, and this is consistent with the irregular flagellum beating reported in~Figure~\ref{fig:helix_3p}. In all three coupling cases, the normalized absolute difference between the actual curvature $\Omega^*$ and the preferred curvature $\Omega$~\eqref{eq:curva} is less than $1\times 10^{-2}$.
%%%

%%%%%%%%%%%%%%%%%%%%%%%%%%%%%%%%%%%%%
\subsubsection{Quasi-planar wave.}
\begin{figure}[!t]
\centering
\includegraphics[width=0.92\linewidth]{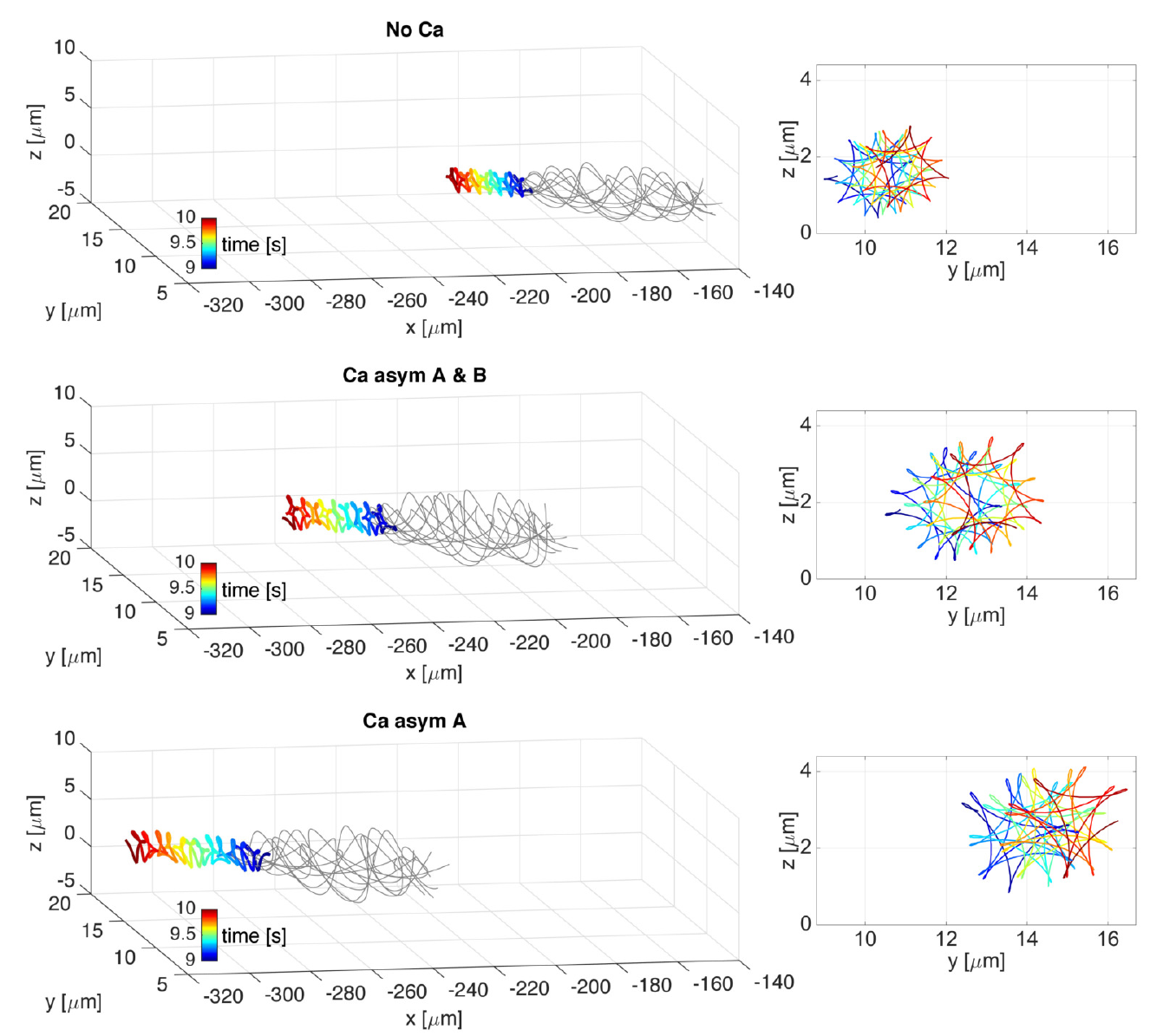}
\caption{\textit{Quasi-planar wave.} 3D overlap of flagellar wave (thin gray lines) and \revtwo{the trajectory of the first point in the direction of swimming} (thick colored line) over the time interval of 1s, form 9s (blue) to 10s (red), on the left. Corresponding flagelloid curve traced by the \revtwo{first point} on the $yz$-plane, on the right. No-coupling case (No Ca, top), asymmetric coupling case (Ca asym $A$ \& $B$, center) and asymmetric coupling only A case (Ca asym $A$, bottom).}
\label{fig:qp_3p}
\vspace*{-9pt}
\end{figure}

In the case of a quasi-planar wave propagating along the flagellum, Figure~\ref{fig:qp_3p} shows flagellar configurations (thin gray lines) and \revtwo{the trajectory of the first point in the direction of swimming} (thick colored line) over the time interval of $1$s, form $9$s (blue) to $10$s (red). 
The three different coupling cases considered are: \textit{No Ca}, \textit{Ca asym $A$ \& $B$} and \textit{Ca asym $A$}. On the right is depicted the corresponding f-curve traced by the \revtwo{rod first point} on the $yz$-plane. In all of the coupling cases, the flagellum configuration shows an irregular beat form, due to the baseline geometric non-linearity introduced by the quasi-planer wave, i.e $A_0\neq B_0$. At the same time, the three cases differ from each other in terms of the emergent shapes of the \revtwo{first point} f-curve, achieved wave amplitude, and linear velocity. In the \textit{No Ca} and \textit{Ca asym $A$ \& $B$} cases, the \revtwo{first point} f-curves {resemble a hypotrochoid with four singular points},
 while in the \textit{Ca asym $A$} case, the \revtwo{first point} f-curves {resemble a hypotrochoid with three singular points}.
Moreover, the two asymmetric coupling cases show visible increased amplitudes in comparison to the no-coupling case. The linear velocity of the quasi-planar wave cases are reported in Table~\ref{tab:vel}. The \textit{Ca asym $A$ \& $B$} and \textit{Ca asym $A$} couplings produce an increase in the linear velocity of $37\%$ and $76\%$
compared to the no-coupling case, respectively.
%%%

\begin{figure}[!t]
\centering
\includegraphics[width=0.8\linewidth]{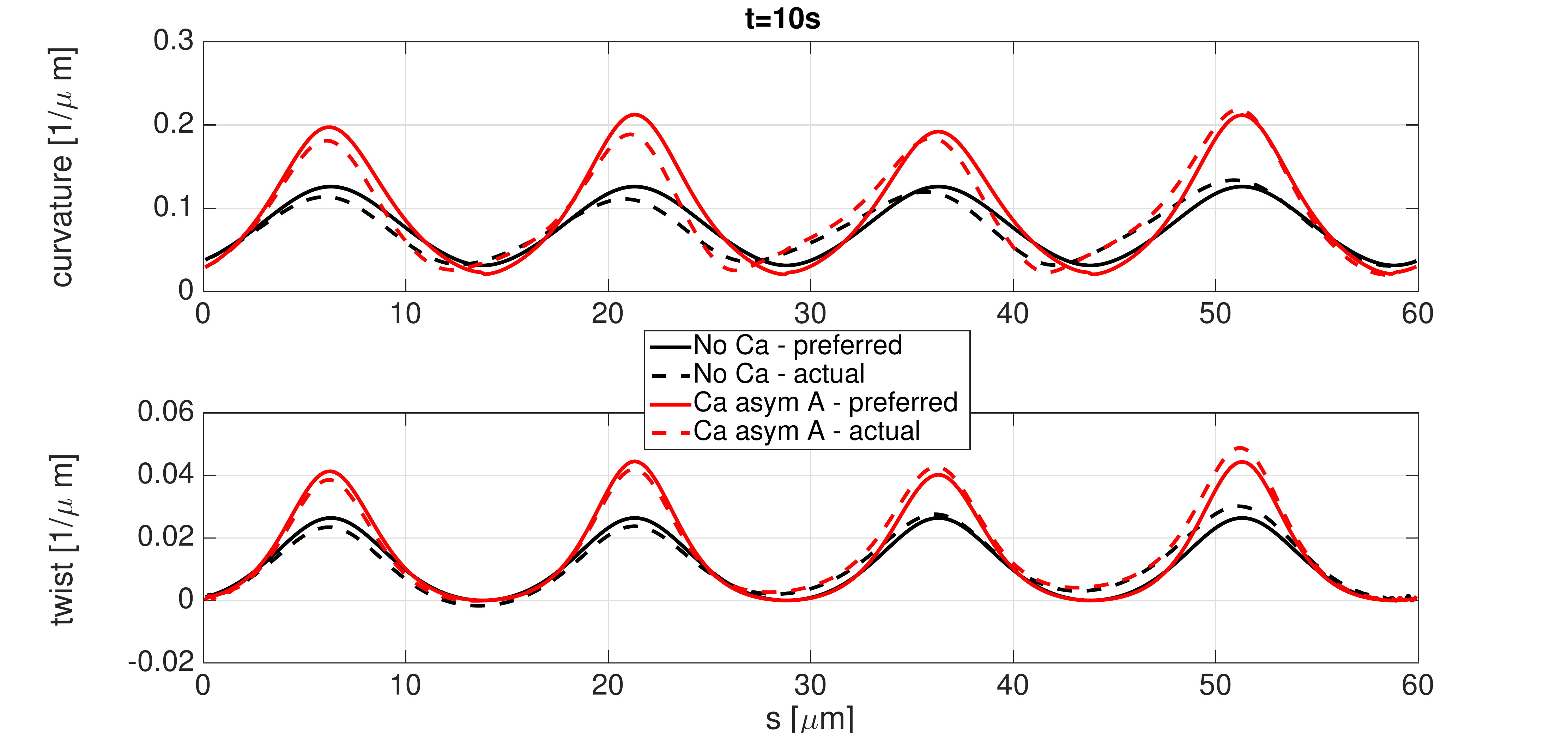}
\caption{\textit{Quasi-planar wave.} Comparison at $t=10$s between preferred (solid lines) and actual (dashed lines) configurations of the rod: curvature in the top panel and twist in the bottom panel for the case of no-coupling (No Ca, black) and in the case of asymmetric coupling only $A$ (Ca asym $A$, red).}
\label{fig:qp_curva}
\vspace*{-9pt}
\end{figure}

In the quasi-planar case, for all three coupling cases considered in Figure~\ref{fig:qp_3p}, the curvature is not constant along the flagellum length and periodically oscillates in time. In particular, Figure~\ref{fig:qp_curva} shows the comparison between preferred curvature $\Omega$~\eqref{eq:curva} and preferred twist $\Omega_3$ (solid lines) with the actual curvature ${\Omega}^*$~\eqref{eq:act_curva} and actual twist $\Omega_3^*$ (dashed lines), for the \textit{No Ca} (black) and the \textit{Ca asym $A$} (red) cases, at time $t=10s$. 
The maximum curvature and twist in the \textit{Ca asym $A$} is almost twice that of the no-coupling case, and both cases present a phase lag between the oscillations of the preferred and computed curvature and twist. 
%%%%

%%%%%%%%%%%%%%%%%%%%%%%%%%%%%%%%%%%%%
\subsubsection{Hypotrochoid approximation of \revtwo{rod first point} f-curves.}
\label{sec:hypo_approx}
As briefly mentioned in the previous sections, the \revtwo{rod first point} f-curves presented in Figures~\ref{fig:helix_3p} and~\ref{fig:qp_3p} resemble hypotrochoid curves.
Hypotrochoid curves can be defined as the trajectories of a point $P$ subjected to a movement composed of two circular motions in opposite directions. The hypotrochoid equation in the complex $xy$-plane can be expressed as
\begin{equation}
z(\gamma) = x(\gamma) + i \ y(\gamma) = \widetilde{R}\exp\left(i\gamma\right) + d \exp\left( -i \frac{\omega_2}{\omega_1} \gamma \right),
\label{eq:hypo}
\end{equation}
where $i$ is the imaginary unit, $\gamma$ is the curve parametrization parameter, $\widetilde{R}$ is the radius of the counterclockwise rotation of frequency $1$~rad/s, and $d$ is the radius of the clockwise rotation of frequency $\omega_2/\omega_1$. For more details on hypotrochoid curves and their parametric representations, please refer to the Appendix~\ref{sec:hypo_app}. Note that in the sperm motility framework, $\omega_1$ represents the sperm roll frequency while $\omega_2$ represents the counter-rotation flagellar frequency~\citep{Woolley98}. 
The shape of the hypotrochoid depends on the ratio between the frequencies $\omega_1$ and $\omega_2$, in particular, the number of singular points $n$ can be determined as
\begin{equation}
n = \Frac{\omega_2}{\omega_1} +1.
\label{eq:hypo_n}
\end{equation}

To better quantify the various f-curve shapes reported in Figures~\ref{fig:helix_3p} and~\ref{fig:qp_3p}, we follow the method reported in~\cite{Woolley98} to approximate the f-curve via a hypotrochoid curve, and the results of the approximation are reported in Figure~\ref{fig:hypo} and in Table~\ref{tab:hypo}. More details on the approximation procedure can be found in the Appendix~\ref{sec:hypo_app}.
Figure~\ref{fig:hypo} shows the comparison between the simulation results for one full rotation of the f-curve (colored curve) and the approximated hypotrochoid curve (black curve) in the case of a helical wave (left), and of a quasi-planar wave (right), considering the asymmetric coupling only A case.
In Table~\ref{tab:hypo}, we report the fitted parameter values for the hypotrochoid in the case of helical and quasi-planar waves and, for the three curvature-calcium coupling cases considered in Figures~\ref{fig:helix_3p} and~\ref{fig:qp_3p}. 
%%%%
\begin{table}[!t] 
\bigskip
\begin{center}
{\mbox{\tabcolsep=5pt\begin{tabular}{@{}clccccccc@{}}
%\tblhead{ 
\toprule
\textit{Wave}&\textit{Coupling case}& 
$\widetilde{R}$ [$\mu$m] & $d$ [$\mu$m] & $\omega_1$[rad/s] & $\omega_2$[rad/s]& $n$\\
\toprule
%[-9.5pt]}\\[-9.5pt]
&No Ca & 
$1.091$ & $0.003$ & $29.7$ & $29.7$ &$2.0^\dag$\\
Helical &Ca asym A\&B & 
$1.451$ & $0.019$ & $20.3$ & $215.6$ & $11.6$\\
&Ca asym A & 
$1.272$ & $0.180$ & $28.5$ & $210.5$ & $8.4$\\
\arrayrulecolor{black!30}\midrule
&No Ca & 
$0.796$ & $0.297$ & $ 62.5$ & $169.0$ & $3.7$\\
Quasi-planar&Ca asym A\&B& 
$1.049$ & $0.375$ & $51.6$  & $172.2$ & $4.3$\\
&Ca asym A & 
$0.945$ & $0.510$ & $68.5$ & $145.9$ & $3.1$\\%[-9.5pt]
\arrayrulecolor{black}\bottomrule
\end{tabular}
\label{tab:hypo}
}}
\caption{The fitted hypotrochoid parameter values varying the flagellar waveform and varying the calcium-curvature coupling condition. No-coupling (No Ca), asymmetric coupling (Ca asym $A$\&$B$) and asymmetric coupling only $A$ (Ca asym $A$). $\dag$Value imposed a priori in order to obtain a circle.}
\end{center}
\end{table}
In the two wave configurations studied, both calcium coupling cases considered produce an increase in the radius $\widetilde{R}$ compared to the non-coupling case, while showing different hypotrochoid shapes. In particular, in the \textit{Ca Asym A\&B} case, $\widetilde{R}$ is about $30\%$ higher than the no-coupling case while maintaining a similar hypotrochoid shape: a circular shape for the helical wave case ($d\simeq0$), and a square shape ($n\simeq4$) for the quasi-planar wave case. However, in the \textit{Ca Asym A} case, compared to the no-coupling case, $\widetilde{R}$ is only about $15\%$ higher, but the hypotrochoid shape varies from a circular ($n=2$) to an octagonal shape ($n\simeq8$) in the helical wave case, and from a square ($n\simeq4$) to a triangular shape ($n\simeq3$) in the quasi-planar wave case. The shapes of the hypotrochoid f-curves extracted from the model results are in agreement with the experimental measurements reported in~\cite{Woolley98,Woolley03} and~\cite{Woolley01}, and reproduced in Figure \ref{exppic}(a). 
In the model results, \revtwo{the left-handedness of the rod trajectories, given by the rod preferred reference configuration~\eqref{eq:wave}, causes a predominant counterclockwise sperm rolling motion, i.e. a predominant counterclockwise rotation of the f-curves reported in Figures~\ref{fig:helix_3p} and~\ref{fig:qp_3p}. Moreover,}
the counter-rotation \revtwo{(clockwise)} of the flagellum is significantly present only if a difference between the preferred amplitudes $A$ and $B$ is introduced. In our model, $A\neq B$ can be obtained in two ways: a baseline geometrical difference in the quasi-planar wave case where $A_0\neq B_0$, and a calcium coupling difference in the \textit{Ca Asym A} case, where only $A$ varies with calcium and $B$ is kept constant to its baseline value. 
%%%%
\begin{figure}[!t]
\centering
\includegraphics[width=0.9\linewidth]{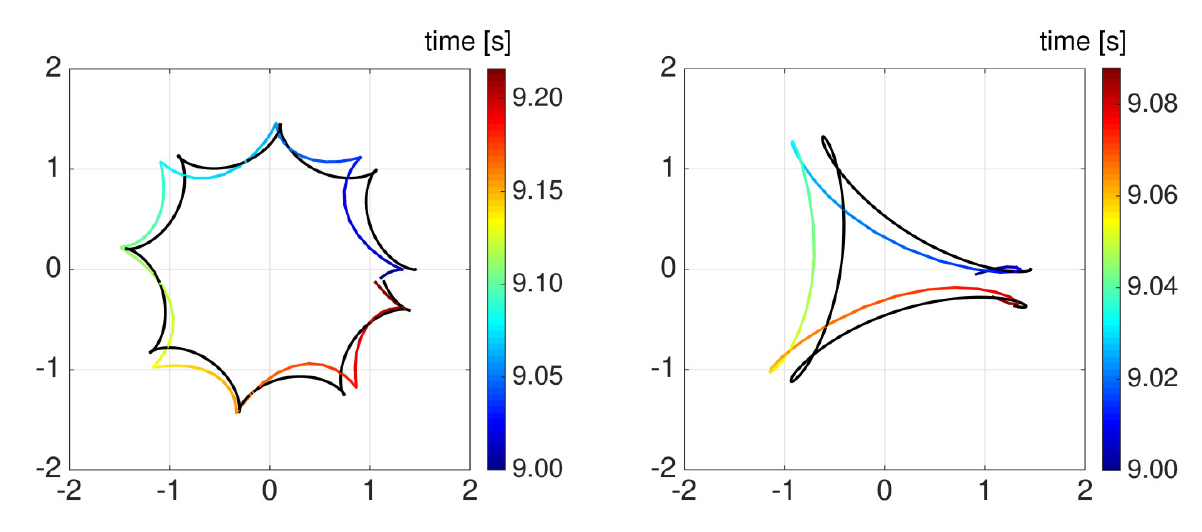}
\caption{\textit{Helical vs Quasi-planar waves.} Comparison between simulation results for a full rotation of the \revtwo{rod first point} f-curve around the center of mass (colored line) and the corresponding fitted hypotrochoid curve (black line) for the case of a helical wave (left) and a quasi-planar wave (right), considering the asymmetric coupling only A case.}
\label{fig:hypo}
\vspace*{-9pt}
\end{figure}

%%%%%%%
\subsubsection{Comparison of helical \& quasi-planar waves.}
\begin{figure}[!t]
\centering
\includegraphics[width=0.95\linewidth]{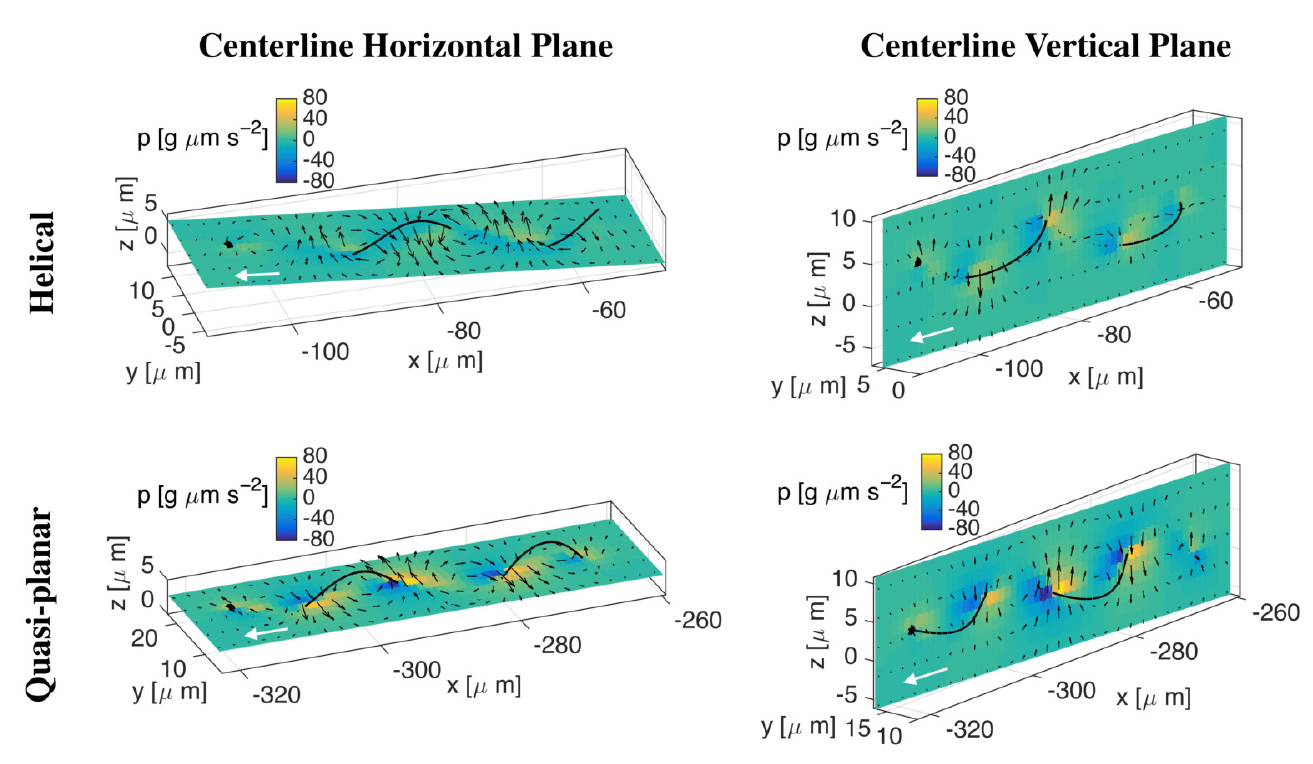}
\caption{\textit{Helical vs Quasi-planar waves.} Fluid velocity fields (black arrows) and pressure ($p$) distributions at time $t=10$s in the flagellum horizontal (left column) or vertical (right column) centerline planes in the case of a helical wave (top row) and a quasi-planar wave (bottom row), for the asymmetric coupling only $A$ case. A filled in sphere is used to denote the \revtwo{rod first point} and \revtwo{a white arrow to denote the} swimming direction.}
\label{fig:3d_vel_p}
\vspace*{-9pt}
\end{figure}

To investigate the fluid-structure interaction between the flagellum and the surrounding fluid, in Figure~\ref{fig:3d_vel_p} we report the fluid velocity fields (black arrows) and pressure ($p$) distributions at time $t=10$s. We choose either a horizontal or vertical plane containing the flagellum centerline in the case of helical and quasi-planar waves, accounting for the \textit{Ca asym $A$} coupling. 
The horizontal centerline plane is defined as the plane passing through the sperm centerline with normal orthogonal to the $y$-axis, while the vertical centerline plane is defined as the plane passing through the sperm centerline with normal orthogonal to the $z$-axis. The pressure distribution range is approximately $-40$ to $40$ g$\mu$ms$^{-2}$ in the helical wave case, and from $-80$ to $80$ g$\mu$ms$^{-2}$ in the quasi-planar wave case. In both wave cases, the velocity fields show vortices along the flagellum length, where the \revtwo{fluid} rotates from regions of negative pressure to regions of positive pressure.
%%%%

In Table~\ref{tab:vel}, we report the values of the linear velocity \revtwo{of the rod first point in the direction of swimming}, rod maximum actual curvature $\Omega^*$, and rod maximum distance from the centerline over the $1$s interval form $9$s to $10$s in the case of no-coupling, asymmetric coupling and asymmetric coupling only $A$, for the helical and quasi-planar waves. 
\begin{table}[!t] 
\bigskip
\begin{center}
{\mbox{\tabcolsep=5pt\begin{tabular}{@{}lcccccc@{}}
\toprule
%\tblhead{ 
\multirow{2}{*}{\makecell[c]{\textit{Coupling case}}}&  \multicolumn{2}{c}{\textit{Linear Velocity} [$\mu$ms$^{-1}$]} &  \multicolumn{2}{c}{\textit{Max. Curvature} [1/$\mu$m]}
&  \multicolumn{2}{c}{\textit{Max. Distance} [$\mu$m]}\\
	& Helical & Quasi-planar & Helical & Quasi-planar & Helical & Quasi-planar\\
	%[-9.5pt]
%	}\\[-9.5pt]
\toprule
No Ca & $8.6$& $22.1$ & $0.09$ & $0.13$ & $2.86$ & $2.87$\\
Ca asym $A$ \& $B$ & $7.6$& $30.4$ & $0.11$ & $0.20$ & $3.94$ & $3.91$\\
Ca asym $A$ & $13.2$ & $38.9$ & $0.16$& $0.22$& $3.82$ &$3.95$\\
\bottomrule
\end{tabular}
\label{tab:vel}
}}
\caption{Comparison of \revtwo{rod first point} linear velocity, maximum curvature, and maximum distance form the centerline for the various flagellar wave and calcium-curvature coupling conditions. No-coupling (No Ca), asymmetric coupling (Ca asym $A$\&$B$) and asymmetric coupling only $A$ (Ca asym $A$).}
\end{center}
\end{table}
In both wave configurations, the  \textit{Ca asym $A$} coupling produces the highest velocities and highest curvatures. 
Note that the linear velocities and maximum curvatures in the quasi-planar wave case are significantly higher than the helical case: from $2.6$ to $4$ times higher for velocities, and from $1.3$ to $1.8$ times higher for curvatures. The maximum distance shows a non-linear behavior as the waveform changes and as the calcium-curvature coupling condition changes. The maximum distance increases approximately $35\%$ in the \textit{Ca asym $A$ \& $B$} coupling case compared to the no-coupling case, and there is no significant variation in the \textit{Ca asym $A$} coupling case compared to the \textit{Ca asym $A$ \& $B$} coupling case. There is no significant difference in the maximum distance values between the helical and the quasi-planar cases, independently from the calcium-curvature coupling considered. \revtwo{Most previous models have indirectly accounted for calcium via a \textit{prescribed} flagellar beat form and hence, we can not compare our emergent amplitudes to these cases. In comparison to models of 2D flagellar bending where calcium dynamics were directly coupled, a non-monotonic behavior in achieved amplitude has also been observed and shown qualitatively in graphs of the flagellar beat form \citep{olson_11,simons_14}. }

\begin{figure}[!t]
\centering
\includegraphics[width=1\linewidth]{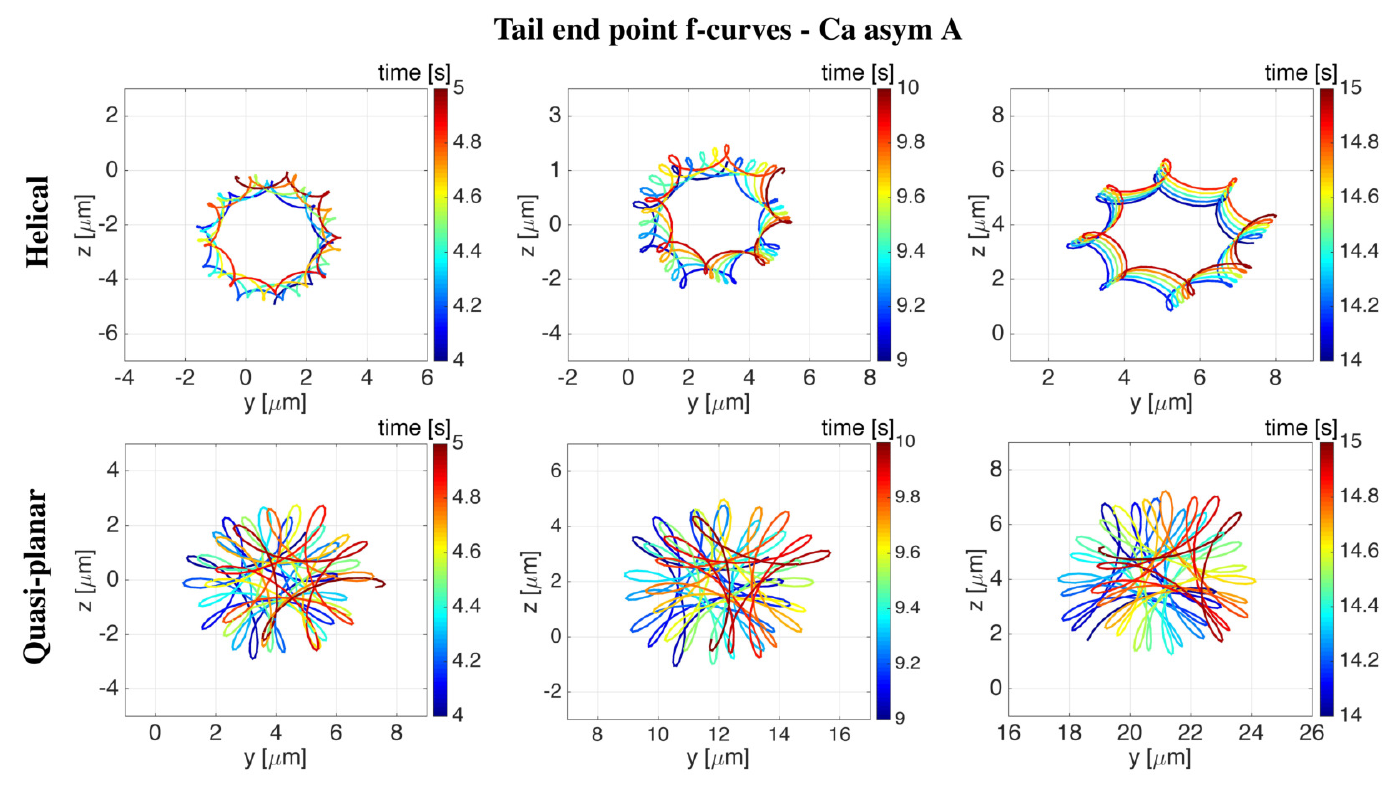}
\caption{\textit{Helical vs Quasi-planar waves.} Tail \revtwo{end point} flagelloid curves (f-curves) on the $yz$-plane over three time intervals of $1$s (left column - from $4$s to $5$s, middle column - from $9$s to $10$s, right column - from $14$s to $15$s) in the case of a helical wave (top row) and a quasi-planar wave (bottom row), for the asymmetric coupling only $A$ (Ca asym $A$).}
\label{fig:trace_tails}
\vspace*{-9pt}
\end{figure}

In Figure~\ref{fig:trace_tails}, the f-curves traced by the end point of the tail on the $yz$-plane are reported over three $1$s time intervals: from $4$s to $5$s (left column), from $9$s to $10$s (middle column) and from $14$s to $15$s (right column).
The helical wave case is reported on the top row and the quasi-planar wave case is on the bottom row, accounting for the \textit{Ca asym $A$} coupling. We observe that the tail f-curves show a similar shape to the f-curves reported in Figures~\ref{fig:helix_3p} and~\ref{fig:qp_3p}, namely {hypotrochoids with eight singular points} 
for the helical wave - \textit{Ca asym $A$} case, and {hypotrochoids with three singular points} 
for the quasi-planar wave - \textit{Ca asym $A$} case. Hence, the results show that the f-curve shape travels along the flagellum length and is conserved in the tail f-curve.
%%%%%

%%%%%%%%%%%%%%%%%%%%%%%%%%%%%
\section{Conclusions}\label{sec:conclusions}
This paper presents the first mathematical model that couples the 3D dynamics of the sperm flagellum and surrounding fluid, with the calcium dynamics inside the flagellum. The coupling is achieved by assuming that the preferred curvature and twist of the flagellum depends on the local calcium concentration. We compare the model results for 2D and 3D flagellar reference beat forms: planar, helical and quasi-planar. Linear velocities, forces and trajectories extracted from the model results are in agreement with experimental data and previously developed mathematical models~\citep{olson_11,simons_14,smith_09,Woolley98,Woolley03,Woolley01}. 

In particular, the planar swimmer results reveal that \textit{(i)},  turning motion on the time scale of seconds is obtained only when asymmetric calcium coupling is considered and \textit{(ii)}, a significantly higher linear velocity is obtained relative to the helical swimmer, in agreement with~\cite{Chwang_71}. \rev{These model results support the hypothesis that helical motion might be an important strategy for chemotaxis sampling. \revtwo{Chemotaxis is defined as the movement towards a chemical concentration gradient such as an egg protein in marine invertebrate sperm or progesterone in human sperm~\citep{Lishko11,Wood_05}.} Swimming with a helical beat form could potentially help the sperm to slow down and sample a larger area, which might aid in the determination of the direction corresponding to an increase in the chemoattractant~\citep{Guerrero11,Su12,Su13}.} In the case of 3D flagellar beat forms, if an asymmetry between the amplitudes $A$ and $B$ is introduced (corresponding to a chirality parameter $0<\alpha<1$ for $A=\alpha B$), the flagelloid curves (f-curves) extracted from the model resemble hypotrochoid curves. In this paper we investigate two sources of asymmetry: a geometrical asymmetry in the quasi-planar wave case, and a calcium coupling asymmetry in the \textit{Ca Asym A} case. The shape of the hypotrochoid f-curve, i.e. the number $n$ and the presence or absence of cusps and self intersections, varies with the choice of preferred flagellar amplitudes $A_0$ and $B_0$, and with the calcium coupling considered. However, the shape of the f-curve is maintained along the length of the flagellum in each case considered. We remark that hypotrochoid f-curves similar to the one extracted from the model results have been obtained in the experimental measurements reported by~\cite{Woolley98,Woolley03} and theoretical predictions reported by~\cite{Smith_09b} and \cite{Ishimoto18}.

Throughout the development of the model, several assumptions were made. We utilized a simplified sperm representation which neglected the cell body or head. Experiments have shown that headless mammalian sperm can still swim in fluid flows~\citep{Miki16}. \rev{Moreover, since the head is small compared to the length scale of the flagellum, a dimensional analysis by \cite{Ishimoto_15} showed that neglecting the head still results in an angular and linear velocity of the correct scale. Therefore, accounting for the head in the model should not significantly affect the overall sperm motility. Further studies are required to investigate the head effect on sperm trajectories and f-curves since the presence of the cell body will certainly increase the drag and potentially bias trajectories and/or turning behavior. In addition, we have used Stokes equations for the fluid flow and all results can be interpreted for the case of a homogeneous fluid with the same viscosity of water. Hence, we have ignored non-Newtonian or resistive effects due to polymers or proteins that are often found in oviductal or vaginal fluid \citep{Miki16,Rutllant05}, as well as polyacrylamide or methylcellulose gels used in experiments \citep{smith_09,Suarez92,Woolley03}. We note that experimental observations where the flagellar beat form varies from 3D to 2D as viscosity varies \citep{smith_09,Woolley03} are not captured in this model, but this will be the focus of future modeling studies using more complicated fluid models. }

\revtwo{In this work, we have focused on constant beat form parameters and material parameters (moduli $a_i$ and $b_i$) in a physiologically representative range (see references in Table \ref{tab:param}). The emergent beat form and trajectory are achieved due to a combination of the preferred beat form parameters, material parameters, as well as contributions from the viscous fluid (via the no-slip condition).
For a preferred symmetric beat form, we observe linear trajectories where there is a slight tilt or yaw in the upward direction (increasing in $y$ for the planar case, Figure \ref{fig:planar_head}, and increasing in both $y$ and $z$ for the quasi-planar and helical cases, Figures \ref{fig:helix_3p} and \ref{fig:qp_3p}), similar to previous models \citep{gillies_09,olson_11,simons_14,Smith_09b}. We note that the emergent flagellar dynamics might exhibit behavior deviating far from the preferred curvature (e.g. buckling instability as in \cite{lee_14} and \cite{Park17}) for a different range of parameters or in the case of anisotropies in material or geometric parameters. To simplify the representation of the flagellum, we have considered an isotropic rod, and neglected the fact that the sperm flagellum has a non-constant cross-sectional radius and  non-constant material properties~\citep{Lesich08,Schmitz04}. Accounting for these anisotropic effects in the mathematical framework presented here introduces interesting challenges from the numerical and modeling view point. Assuming a non-constant rod radius implies that the regularization parameter $\epsilon$ decreases along the rod length. If $\epsilon$ is too small compared to the spatial discretization step $\triangle s$, numerical instabilities might occur, i.e. the fluid might leak through the rod centerline. In order to prevent these instabilities, we could either decrease $\triangle s$ and at the same time increase the computational costs, or we could adopt a more stable regularized Stokeslet framework recently presented in~\cite{cortez_18}.
In the same fashion of~\cite{olson_11}, we could account for the material anisotropy by including a taper in the material properties from the base of the flagellum to the tip of the tail. However, in~\cite{olson_11}, a 2D Euler elastica model was used where only two material parameters were included, i.e. bending and tensile stiffness; the bending stiffness was fixed, while a linear taper or a 4th order taper was used for the tensile stiffness. In this work we utilize a 3D Kirchhoff rod model that includes six material parameters, $a_i$ and $b_i$ for $i=1,2,3$, and careful consideration and possibly additional experimental measurements would need to be made to account for the taper effect in each of the moduli. It would be interesting to investigate these anisotropic effects coupled with calcium dynamics, since they might predominantly affect sperm trajectories in the end piece region where the flagellum is thinner \citep{cummins_85}.
}

\revtwo{In contrast to previous models that have indirectly coupled calcium to} hyperactivated motility\revtwo{,  we have directly accounted for this in mammalian sperm via} an evolving calcium concentration that is mediated by CatSper channels, ATP-ase pumps, and a calcium store (redundant nuclear envelope) where calcium release is coupled to the evolving concentration of inositol 1,4,5-trisphosphate or $IP_3$. \revtwo{We note that this modeling framework could be utilized to investigate how different pharmaceutical treatments could enhance or diminish mammalian sperm motility by varying parameters related to the different calcium pumps and channels \citep{Lishko12,Shukla12,Zheng13}.}
However, in other species, different calcium signaling mechanisms might play an important role. For example, in sea urchin sperm, an increase in internal calcium is associated with an increase in internal cGMP (cyclic guanosine monophosphate) and the calcium channels could be CatSper or other voltage sensitive calcium channels~\citep{Darson_08,Seifert15,Wood_05}. Hence, in order to study the effect of calcium dynamics on sperm motility in different species, the calcium model used here should be adapted to account for the species specific calcium signaling pathways and associated channels. We remark that, in addition to the CatSper signaling pathways related to hyperactivation in mammalian sperm, there are other mechanisms that could guide the sperm to the egg and change the flagellar waveform. These include chemotaxis, i.e. the effect of a chemical concentration gradient~\citep{Lishko11,Wood_05}, rheotaxis, i.e. the effect of a background flow~\citep{Miki16}, and thermotaxis, i.e. the effect of the fluid temperature gradient~\citep{Boryshpolets15}. It would be an interesting research direction to expand the model in the current study to include the effect of a background flow and an evolving gradient of chemoattractants or calcium in the fluid. 

%%%%%%%%%%%%%%%%%%%%%%%%%%%%%
\section*{Acknowledgments}
The work of L. Carichino and S.D. Olson was supported, in part, by NSF DMS 1455270. All simulations were run on a cluster acquired through NSF DMS 1337943. \rev{We also wish to thank Padraig \'{O} Cath\'{a}in and Jianjun Huang for many useful discussions.}

%%%%%%%%%%%%%%%%%%%%%%%%%%%%%
\newpage
\begin{appendices}
\section*{Appendices}
\addcontentsline{toc}{section}{Appendices}
\renewcommand{\thesubsection}{\Alph{subsection}}
\renewcommand{\theequation}{\thesubsection.\arabic{equation}}
\vspace{-0.5cm}
\rev{\subsection{Rod reference configuration}\label{AppCurva}
Given the preferred reference configuration of $\widehat{\vector{X}}(t,s)$ in~\eqref{eq:wave} with components $x(t,s)$, $y(t,s)$, and $z(t,s)$ for arc length parameter $s$ at a given time $t$, we can determine equally spaced points along the curve. These points are then used to define the corresponding triad $\{\widehat{\vector{D}}^1(t,s), \widehat{\vector{D}}^2 (t,s), \widehat{\vector{D}}^3 (t,s) \}$. According to the arc length formula, we have 
\begin{equation}\label{eq:arc_length}
s = \dint_0^{s} \sqrt{(\partial_\xi x(t,\xi))^2 + (\partial_\xi y(t,\xi))^2 + (\partial_\xi z(t,\xi))^2} d\xi.
\end{equation}
Applying the chain rule to $y(t,s)$ and $z(t,s)$, we have 
\begin{equation}\label{eq: yzprime}
\partial_s y(t,s)= A k \cos(k x(t,s) -\sigma t)\cdot \partial_s x(t,s) \quad \mbox{ and } \quad \partial_s z(t,s)= -B k \sin (k x(t,s) -\sigma t)\cdot \partial_s x(t,s).
\end{equation}
Differentiating the arc length formula in~\eqref{eq:arc_length} with respect to $s$ and plugging the formulae for $\partial_s y$ and $\partial_s z$ from~\eqref{eq: yzprime} into the equation, we can then solve for $\partial_s x$ as
\begin{equation}
\partial_s x(t,s) = \frac{1}{\sqrt{1 + A^2 k^2 \cos^2(k x(t,s) -\sigma t) +  B^2 k^2 \sin^2(k x(t,s) -\sigma t) }}.
\label{equ: xprime}
\end{equation}
The discretization of the rod in terms of $x(t,s)$ (and hence $y(t,s)$ and $z(t,s)$) can then be obtained by solving this differential equation. The tangent vector is $\widehat{\vector{D}}^3$ and is initialized as $\widehat{\vector{D}}^3= (\partial_s x,\; \partial_s y,\; \partial_s z)$.  \revtwo{The normal direction $\widehat{\vector{D}}^1=(\vector{e}_3 \times \widehat{\vector{D}}^3)/(|\vector{e}_3 \times \widehat{\vector{D}}^3|)$ where $\vector{e}_3= (0,0,1)$, which means that $\widehat{\vector{D}}^1 = \Frac{(-\partial_s y,\; \partial_s x,\; 0)}{\sqrt{(\partial_s x)^2 + (\partial_s y)^2}}$. 
We remark that given the form of the preferred reference configuration adopted in~\eqref{eq:wave}, the tangent vector $\widehat{\vector{D}}^3$ is never parallel to $\vector{e}_3$, i.e. $\partial_s x \neq 0 $ and $\partial_s y \neq 0$, and this guarantees that $\widehat{\vector{D}}^1$ is well defined.}
The binormal is then $\widehat{\vector{D}}^2=\widehat{\vector{D}}^3 \times \widehat{\vector{D}}^1$.
Once the orthonormal triad is established as given, we can obtain the preferred strain twist vector $(\Omega_1,\Omega_2,\Omega_3)$ using $\Omega_i=\partial_s\widehat{\vector{D}}^j\cdot\widehat{\vector{D}}^k$ for $(i,j,k)$ a cyclic permutation of (1,2,3).}
%%%%%%%%%%%%%%%%%%%%%%%%%%%%%
\subsection{Hypotrochoid approximation}\label{sec:hypo_app}
\begin{figure}[!t]
\centering
\includegraphics[width=0.4\linewidth]{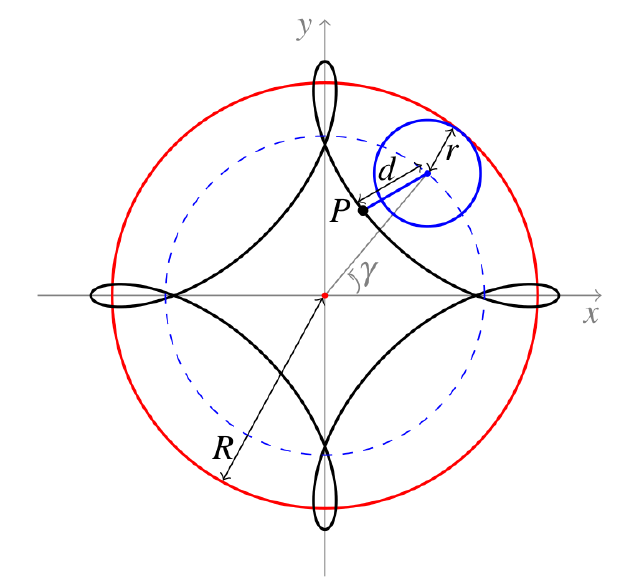}
\caption{Sketch of a theoretical hypotrochoid curve (black) obtained by tracing a point $P$ attached to a small circle (blue) of radius $r$ rolling around the inside of a bigger circle (red) of radius $R$, where $d$ is the distance between the point $P$ and the center of the small circle. $\gamma$ is the angle formed by the horizontal axis and the center of the small rolling circle. The blue dashed line represents the trace of the center of the small circle, i.e. a circle of radius $\widetilde{R}=R-r$. $R=2$, $r=0.5$ and $d=0.7$.}
\label{fig:appendix}
\vspace*{-9pt}
\end{figure}

Hypotrochoid curves are usually defined as the curve created by tracing a point $P$ rigidly attached to a small circle of radius $r$ rolling around the inside of a bigger circle of radius $R$, where $d$ is the distance between the point $P$ and the center of the small circle, as sketched in Figure~\ref{fig:appendix}. 
The parametric equations for a hypotrochoid in the $xy$-plane can be expressed as
\begin{equation}
\label{eq:hypo_app}
x(\gamma)= (R-r)\cos \gamma + d \cos\left( \Frac{R-r}{r} \gamma \right), \qquad
y(\gamma)= (R-r)\sin \gamma - d \sin\left( \Frac{R-r}{r} \gamma \right),
\end{equation}
where $\gamma$ is the angle formed by the horizontal axis and the center of the small rolling circle. 
Equation~\eqref{eq:hypo_app} can be derived directly from~\eqref{eq:hypo}, given 
\begin{equation}
\widetilde{R}=R-r \qquad \mbox{ and } \qquad n = \Frac{R}{r} = \Frac{\omega_2}{\omega_1} +1.
\end{equation}
Hence, hypotrochoid curves can also be described as the trajectories of a point $P$ subjected to a movement composed of two circular motions in opposite directions: the center of the small circle is rotating around the origin with a counterclockwise circular motion of radius $\widetilde{R}=R-r$ and frequency $1$rad/s (represented by the blue dashed line in Figure~\ref{fig:appendix}), and the point $P$, rigidly attached to the small circle, is rotating around the center of the small circle with a clockwise circular motion of radius $d$ and frequency $\omega_2/\omega_1$. As mentioned in Section~\ref{sec:hypo_approx}, this last characterization of hypotrochoid curves can be directly interpreted in the sperm motility framework. Note that the curve corresponding to~\eqref{eq:hypo} is equivalent to the curve generated by a counterclockwise circular motion of radius $\widetilde{R}$ and frequency $\omega_1$ together with a clockwise circular motion of radius $d$ and frequency $\omega_2$; $P$ is just going across them with a different speed as the parameter $\gamma$ varies.
%%%%

Fixing $n$, the value of $d$ determines if singular points are present, i.e. cusps or self intersections (crunodes). We remark that in Section~\ref{sec:hypo_approx} we refer to $n$ as the number of singular points of the curve, however this is a slight abuse of notation since, fixing $n$, for some values of $d$ the curve does not present any singular point and can be described as a rounded approximation of a $n$-sided regular polygon.
Note that $d$ admits positive and negatives values: at  $\gamma=0$, if $d>0$  the point $P$ is chosen to the right of the small circle center, otherwise if $d<0$ the point $P$ is chosen to the left. Moreover, the point $P$ can be chosen to be either inside ($|d|<r$), outside ($|d|>r$) or on the circumference of the small circle ($|d|=r$). The case of $d=r>0$ corresponds to hypocycloid curves, while the limit case of $d=0$ corresponds to a circle of radius $\widetilde{R}$.

Fitting sperm trajectory data to hypotrochoid curves involves different mathematical challenges. Starting from finding the best minimization algorithm and initial value to approximate the three parameters, $(\widetilde{R},d,n)$ or $(R,r,d)$, that uniquely identify the curve. To our knowledge, there is only one result available in the literature on least-squares hypotrochoid curve fitting in \cite{Sinnreich_16}, where a method for a particular case of hypotrochoid curve is presented, i.e. rounded approximation of regular polygons with no cusps and no self-intersections. A further challenge comes form the fact that for each data point the corresponding $\gamma$ is also an unknown, since by definition $\gamma$ is in general not equal to the point polar angle $\theta$. We can relate $\theta$ and $\gamma$ using the inverse tangent of the ratio between the $y$ and $x$ coordinates of the data point. However, this correspondence is not one-to-one when self intersections are present in the hypotrochoid curve. For this reason, we chose to use the method reported in~\cite{Woolley98} to approximate the simulation f-curves with a hypotrochoid curve. This method exploits the definition of a hypotrochoid curve as two circular motions, and uses geometric and physical principles to estimate the curve parameters without utilizing a minimization algorithm, and can be used for approximating any kind of hypotrochoid curves, including ones that exhibit self-intersections. 

Given the f-curves data points, we follow~\cite{Woolley98} to find the approximating hypotrochoid curve as detailed below.
\begin{enumerate}
\item Starting at time $t=9s$, the roll frequency $\omega_1$ is estimated as the frequency of a full rotation of the \revtwo{rod first point} around its center of mass;
\item $n$ is estimated as the ratio between $2\pi$ and the mean angular separation between two singular points over a full rotation;
\item The flagellar frequency is $\omega_2 = \omega_1(n-1)$;
\item Let $d_{min}$ and $d_{max}$ be respectively the average minimum and maximum distance from the center of mass over one rotation, then the remaining hypotrochoid parameters can be estimated as follows
\begin{equation}
d = \Frac{d_{max} - d_{min}}{2}, \qquad R = \Frac{\omega_1 + \omega_2}{\omega_2} \left( d_{min} + d\right), \qquad r=\Frac{\omega_1}{\omega_2} \left( d_{min} + d\right).
\end{equation}
\end{enumerate}
\end{appendices}

%%%%%%%%%%%%%%%%%bibliography style

%%%%%%%%%%%%%%%%%%%%%%%
%%%%%%%%%%%%%%%%%%%%%%%%%%%%%%%%%%

\end{document}